\documentclass[9pt,11pt]{article}%
\usepackage{amsfonts}
\usepackage[applemac]{inputenc}
\usepackage{amssymb}
\usepackage{makeidx}
\usepackage{graphicx}
\usepackage{amsmath}
\usepackage[a4paper]{geometry}
\usepackage[displaymath]{lineno}
\usepackage[singlespacing]{setspace}
\usepackage[normalem]{ulem}
\usepackage{lineno}
\usepackage{color}%
\providecommand{\U}[1]{\protect\rule{.1in}{.1in}}
\providecommand{\U}[1]{\protect\rule{.1in}{.1in}}

\let\epsilon\varepsilon

\newtheorem{thm}{Theorem}
\newtheorem{lem}[thm]{Lemma}

\newtheorem{cor}[thm]{Corollary}
\newtheorem{rem}{Remark}

\newtheorem{defi}{Definition}

\newenvironment{dem}[1][Proof]{\noindent \textbf{#1.} }{\ \rule{0.5em}{0.5em}}
\begin{document}

\title{A general representation\ of $\delta$-normal sets to sublevels of convex
functions\thanks{This work is partially supported by CONICYT grant Fondecyt
1151003, Conicyt-Redes no. 150040, and Mathamsud 17-MATH-06.}}
\author{Abderrahim Hantoute\thanks{e-mail: ahantoute@dim.uchile.cl}, Anton
Svensson\thanks{e-mail: asvensson@dim.uchile.cl}\\{\scriptsize {Center for Mathematical Modeling, Department of Mathematical
Engineering,}}\\{\scriptsize Universidad de Chile, Santiago, Chile} \\\quad\\\quad\\Dedicated to Prof. Michel Théra on his 70th birthday}
\date{}
\maketitle

\begin{abstract}
The $(\delta$-) normal cone to an arbitrary intersection of sublevel sets of
proper, lower semicontinuous, and convex functions is characterized, using
either $\varepsilon$-subdifferentials at the nominal point or exact
subdifferentials at nearby points. Our tools include ($\varepsilon$-) calculus
rules for sup/max functions. The framework of this work is that of a locally
convex space, however, formulas using exact subdifferentials require some
restriction either on the space (e.g. Banach), or on the function (e.g. epi-pointed).



\end{abstract}

\section{Introduction}

In \cite{cabot2014sequential}\textbf{ }the authors prove that for a convex,
proper and lower semicontinuous (lsc) function $\Phi:X\rightarrow
\mathbb{R}\cup\left\{  +\infty\right\}  ,$ defined on a (reflexive) Banach
space $X,$ the normal cone to the sublevel set $[\Phi\leq\Phi(\bar
{x})]:=\{x\in X\mid\Phi(x)\leq\Phi(\bar{x})\}$ at a point $\bar{x}$ in the
effective domain of $\Phi,$ is characterized as the (norm) upper-limit of
directions in the Fenchel subdifferential of $\Phi$ at sufficiently close
points to $\bar{x};$ that is, \
\begin{equation}
\mathrm{N}_{[\Phi\leq\Phi(\bar{x})]}(\bar{x})=\limsup_{\substack{x\rightarrow
\bar{x}}}\mathbb{R}_{+}\partial\Phi(x).\label{thibault}%
\end{equation}
This relation is also valid in general Banach spaces but up to some
specifications of both the $\limsup$ and the convergence of $x$ to $\bar{x}.$
The proof of (\ref{thibault}) and its generalization to Banach spaces given in
\cite{cabot2014sequential}, is based on sequential calculus rules for the
subdifferential of composite functions developed in
\cite{thibault1997sequential}. In this work, assuming that $X$ is a general
locally convex (lc) space with a topological dual $Y,$ we propose another
route to approach this problem and provide different\ characterizations
of\ the $\delta$-normal set to an arbitrary intersection of\ sublevel sets
\[
\mathrm{N}_{\cap_{t\in T}[\Phi_{t}\leq\Phi_{t}(\bar{x})]}^{\delta}(\bar
{x}):=\{x^{\ast}\in Y\mid\left\langle x^{\ast},x-\bar{x}\right\rangle
\leq\delta\text{ for all }x\in\cap_{t\in T}[\Phi_{t}\leq\Phi_{t}(\bar
{x})]\},\text{ }\delta\geq0,
\]
for convex functions $\Phi_{t}$ which are indexed in an arbitrary set $T$.
Compared to formula (\ref{thibault}), the present\ characterization only
involves the reference point $\bar{x}$ rather than nearby points, and uses
$\varepsilon$-subdifferentials instead of exact ones.

Owing to Brøndsted-Rockafellar's theorem (\cite{borwein1982note}), these
epsilon-like formulas for $\delta$-normal sets\ easily lead to
characterizations in the line of formula (\ref{thibault}). This passage from
$\varepsilon$-subdifferentials to exact ones will require\ some natural
conditions either on the underlying space or on the associated function;
typically, that either $X$ is Banach or that $\Phi$ satisfies some
continuity/coercivity conditions.

More generally, we shall prove that if $\Phi=\sup_{t\in T}\Phi_{t},$ with
$\Phi_{t}:X\rightarrow\mathbb{R}\cup\left\{  +\infty\right\}  ,$ $t\in T,$
being convex, proper and lsc, then for every $\delta\geq0,$ $\bar{x}%
\in\operatorname*{dom}\Phi$ and $\lambda\in\left(  -\infty,+\infty\right]  $
we have that
\begin{equation}
\mathrm{N}_{([\Phi\leq\lambda]\cap\operatorname*{dom}\Phi)\cup(\bar{x}%
+[\Phi\leq\Phi(\bar{x})]_{\infty})}^{\delta}(\bar{x})=\limsup_{\substack{\sum
\nolimits_{i\in\overline{1,k}}\mu_{i}(\nu-\Phi_{t_{i}}(\bar{x})+\varepsilon
_{i})\rightarrow\delta,\text{ }\nu\uparrow\lambda\text{ }\\\mu_{i},\text{
}\varepsilon_{i}\geq0,\text{ }t_{i}\in T,\text{ }k\in\mathbb{N}}%
}\sum\nolimits_{i\in\overline{1,k}}\mu_{i}\partial_{\varepsilon_{i}}%
\Phi_{t_{i}}(\bar{x}), \label{our}%
\end{equation}
where the subscript $\infty$ refers to the recession cone,
and\ $\operatorname*{dom}\Phi$ denotes the effective domain of $\Phi.$ The
sets $\operatorname*{dom}\Phi$ and $\bar{x}+[\Phi\leq\Phi(\bar{x})]_{\infty}$
are superfluous when $\Phi(\bar{x})\leq\lambda<\infty,$ so that the left-hand
side in (\ref{our}) reduces to the usual $\delta$-normal set of $[\Phi
\leq\lambda].$ However, these two extra sets are necessary to make formula
(\ref{our}) meaningful when $\bar{x}$ may lie outside the set $[\Phi
\leq\lambda];$ for instance, when $[\Phi\leq\lambda]=\emptyset.$ The case
$\lambda=+\infty$ is also meaningful since it leads to a characterization of
the normal cone to the domain of $\Phi.$ It is worth recalling that when
$\Phi_{t}\equiv\Phi$ and $\Phi(\bar{x})\leq\lambda,$ formula (\ref{our})
yields a well-known result, which was first established in \cite{kutateladze}
(see, also, \cite{hiriart1995subdifferential}). Formula (\ref{our}) also
naturally simplifies under Slater's type conditions, giving rise to familiar
results (see \cite{rockafellar2015convex}).

The passage from single to arbitrary intersections of sublevel sets will be
made possible through the investigation in this work of new and general
($\varepsilon$-)subdifferential calculus rules for pointwise suprema,
extending\ some previous results on this theme (\cite{Brondsted72,
CorreaHantLopez, HanLop08, HanLopZal2008, ioffe11, LiNg11, Volle12}). These
rules are extensively used at different stages of the proof of (\ref{our}).

The going back and forth between (\ref{thibault}) and (\ref{our}) will be made
clear through the use of Brøndsted-Rockafellar's like-theorem, which allows
rewriting\ (\ref{our}) by means only of exact subdifferentials. However, this
approach requires some extra conditions, in light of the example of proper,
convex and lsc functions, defined in non-complete normed spaces, which have
an\ empty subdifferential at every point
(\cite{brondsted1965subdifferentiability}). Despite this limitation, we shall
see that the following relation
\[
\mathrm{N}_{[\Phi\leq\Phi(\bar{x})]}(\bar{x})=\limsup_{\substack{x\rightarrow
_{\operatorname*{dom}\Phi^{\ast}}\bar{x},\text{ }\mu\geq0\\\mu(\Phi
(x)-\Phi(\bar{x}))\rightarrow0\\\mu\left\langle \cdot,x-\bar{x}\right\rangle
\rightarrow0}}\mu\partial\Phi(x),
\]
still hold true for a wide class of functions, which includes for example
functions defined on locally convex spaces which either they or their
conjugates are finite and continuous at some point. This class, referred to as
the class of epi-pointed functions\ (e.g., \cite{correa2016bronsted}), also
contains convex functions defined in Banach spaces up to some appropriate localization.

The need for explicit characterizations of the normal cone to sublevel sets is
fundamental in optimization theory, namely in the derivation of optimality
conditions for convex programming problems (e.g.,
\cite{rockafellar2009variational, rockafellar2015convex}). It is also relevant
in the investigation of stationarity and stability properties of different
dynamic systems (\cite{alvarez2000heavy, cabot2009long, cabot2014sequential}).

The previous formulas are applied at the end of this work to
spectral\ functions (\cite{lewis, lewisin}). As it is expected, the associated
normal cone will only depend on the values of the function on the range of
eigenvalues vectors.

The paper is organized as follows: In section \ref{sect1} we give the
necessary definitions and basic notations. In section \ref{sect2} we develop
some approximate subdifferential calculus rules for the max function, which
are needed for our analysis. In section \ref{sect3} we give in Theorem
\ref{tmain} the main characterization of the normal cone to\ sublevel sets,
which uses the approximate subdifferential. This characterization will be
rewritten in section \ref{sect4}, Theorem \ref{rii}, by means only of exact
subdifferentials. Finally, in section \ref{sect5}, we apply the previous
result to investigate the case of spectral functions.

\section{Notations and Definitions\label{sect1}}

We recall in this section some definitions and notations that will be used in
the sequel. We consider a dual pair $(X,Y)$ of real locally convex (lc, for
short) spaces $X$ and $Y,$ defined via a bilinear form $\left\langle x^{\ast
},x\right\rangle :=\left\langle x,x^{\ast}\right\rangle :=x^{\ast}(x)$,
$x^{\ast}\in Y$, $x\in X.$ By $\mathcal{N}_{X}(x)$ we refer to the family of
absolutely convex neighborhoods of $x.$ The origin vectors are denoted by
$\theta.$

Given a non-empty set $S\subset X$, by $\overline{S}$ (or $\operatorname*{cl}%
S),$ and $\operatorname*{aff}S$ we denote the closure, and the affine hull of
$S,$ respectively. The relative interior of $S$ is the interior of $S$
relative to $\operatorname*{aff}S$ when this set is closed, and the emptyset
otherwise. The polar set, the dual cone, and the orthogonal space of $S$ are
the subsets of $Y$ given by $S^{\circ}:=\{x^{\ast}\in Y\mid\left\langle
x^{\ast},x\right\rangle \leq1$ for all $x\in S\}$, $S^{-}:=\{x^{\ast}\in
Y\mid\left\langle x^{\ast},x\right\rangle \leq0$ for all $x\in S\},$ and
$S^{\perp}:=\{x^{\ast}\in Y\mid\left\langle x^{\ast},x\right\rangle =0$ for
all $x\in S\}$, respectively.$\ $Given $\delta\geq0,$ the $\delta$-normal set
to $S$ at $\bar{x}$ is the set
\[
\mathrm{N}_{S}^{\delta}(\bar{x}):=\left\{  x^{\ast}\in Y\ |\ \left\langle
x^{\ast},x-\bar{x}\right\rangle \leq\delta,\ \forall x\in S\right\}  ;
\]
we call $\mathrm{N}_{S}(\bar{x}):=\mathrm{N}_{S}^{0}(\bar{x})$ the normal
cone\ to $S$ at $\bar{x}$.

We fix a function $\Phi:X\rightarrow\mathbb{R}\cup\left\{  +\infty\right\}  .$
We say that $\Phi$\ is proper if its (effective) domain, $\operatorname*{dom}%
\Phi:=\left\{  x\in X\mid\Phi(x)<\infty\right\}  ,$ is nonempty; convex (lower
semi-continuous (lsc), resp.), if its epigraph,
\[
\operatorname*{epi}\Phi:=\left\{  (x,\lambda)\in X\times\mathbb{R}\mid
\Phi(x)\leq\lambda\right\}  ,
\]
is convex (closed, resp.). If $\Phi$ is proper, convex and lsc we write
$\Phi\in\Gamma_{0}(X).$ The restriction function of $\Phi$ to $S$ is denoted
by $\Phi_{\mid S}.$ The closed convex hull of $\Phi$ is defined as%

\[
\overline{\operatorname*{co}}\Phi(x)=\liminf_{u\rightarrow x}\inf\left\{
\sum_{i\in\overline{1,k}}\lambda_{i}\Phi(x_{i})\mid u=\sum_{i\in\overline
{1,k}}\lambda_{i}x_{i},\ (\lambda_{i})\in\Delta_{k},\ k\in\mathbb{N}\right\}
,\text{ }x\in X,
\]
where $\Delta_{k}:=\{(\lambda_{1},\cdots,\lambda_{k})\mid\lambda_{i}>0$,
$\lambda_{1}+\cdots+\lambda_{k}=1\}.$ The sublevel set of $\Phi$ at
$\lambda\in\left]  -\infty,+\infty\right]  $ is the set
\[
\lbrack\Phi\leq\lambda]:=\{x\in X\mid\Phi(x)\leq\lambda\}.
\]

Assume that $\Phi\in\Gamma_{0}(X).$ For\ $\varepsilon\geq0$ the $\varepsilon
$-subdifferential of $\Phi$ at $x\in\operatorname*{dom}\Phi$ is
\[
\partial_{\varepsilon}\Phi(x):=\left\{  x^{\ast}\in Y\ |\ \langle x^{\ast
},u-x\rangle\leq\Phi(u)-\Phi(x)+\varepsilon,\,\forall u\in X\right\}  ;
\]
we write $\partial_{\varepsilon}\Phi(x):=\emptyset$ if $\varepsilon<0$ or if
$x\notin\operatorname*{dom}\Phi$. The subdifferential of $\Phi$ at $x$ is the
set $\partial\Phi(x):=\partial_{0}\Phi(x)$. The directional $\varepsilon
$-derivative of $\Phi$ at $x$ in a direction $v\in X$ is defined as
\[
\Phi_{\varepsilon}^{\prime}(x,v):=\inf_{t>0}\frac{\Phi(x+tv)-\Phi
(x)+\varepsilon}{t};
\]
again, if $\varepsilon=0,$ we just call it directional derivative and write
$\Phi^{\prime}(x,v).$ Equivalently, for $\varepsilon>0$ we have \cite[Theorem
2.4.11]{zalinesku2000convex}
\[
\Phi_{\varepsilon}^{\prime}(x,v)=\sup_{x^{\ast}\in\partial_{\varepsilon}%
\Phi(x)}\left\langle x^{\ast},v\right\rangle .
\]
The asymptotic function of $\Phi$, $\Phi^{\infty}:X\rightarrow\mathbb{R}%
\cup\left\{  +\infty\right\}  ,$ is defined as
\[
\Phi^{\infty}(v):=\sup_{t>0}t^{-1}(\Phi(x+tv)-\Phi(x)).
\]
The conjugate of $\Phi$ is\ the proper, convex and lsc function defined on $Y$
as
\[
\Phi^{\ast}(x^{\ast}):=\sup\{\left\langle x^{\ast},x\right\rangle
-\Phi(x),\text{ }x\in X\}.
\]
The indicator function of $S$, $\mathrm{I}_{S}:X\rightarrow\mathbb{R}%
\cup\left\{  +\infty\right\}  ,$ is defined by
\[
\mathrm{I}_{S}(x)=\left\{
\begin{matrix}
0 & \text{if }x\in S\\
+\infty & ~~\text{if }x\notin S,
\end{matrix}
\right.
\]
while the support function, $\mathrm{\sigma}_{S}:X\rightarrow\mathbb{R}%
\cup\left\{  +\infty\right\}  $ (when $S\subset Y$), is defined as the
conjugate of $\mathrm{I}_{S}.$ We shall frequently use the following relation,
which holds for every function $\varphi:X\rightarrow\mathbb{R}\cup\left\{
+\infty\right\}  $ having a proper conjugate,
\begin{equation}
\varphi^{\ast\ast}:=(\varphi^{\ast})^{\ast}=\overline{\operatorname*{co}%
}\varphi. \label{cob}%
\end{equation}
We know that
\begin{equation}
\Phi^{\infty}=\mathrm{\sigma}_{\operatorname*{dom}\Phi^{\ast}}, \label{in1}%
\end{equation}

\begin{equation}
\lbrack\Phi\leq\Phi(x)]_{\infty}=[\Phi^{\infty}\leq0]=[\mathrm{\sigma
}_{\operatorname*{dom}\Phi^{\ast}}\leq0]=(\operatorname*{dom}\Phi^{\ast})^{-};
\label{in}%
\end{equation}
here, $S_{\infty}$ denotes the asymptotic cone of $S$ defined via the
relation\ $\mathrm{I}_{S_{\infty}}=(\mathrm{I}_{S})^{\infty}.$

Recall that for any set $T$ and $\Phi_{t}\in\Gamma_{0}(X),$ $t\in T,$ we have
that (see \cite{moreau})\
\begin{equation}
(\sup_{t\in T}\Phi_{t})^{\ast}=\overline{\operatorname*{co}}(\inf_{t\in T}%
\Phi_{t}^{\ast}),\text{ }(\inf_{t\in T}\Phi_{t})^{\ast}=\sup_{t\in T}\Phi
_{t}^{\ast}.\label{in3}%
\end{equation}

Finally, given a multifunction $M:U\rightrightarrows V,$ defined between two
topological spaces $U$ and $(V,\tau),$ the Painleve-Kuratowski upper limit of
$M$ at $\bar{u}\in U$ is defined as
\[
\tau\text{-}\limsup\limits_{u\rightarrow\bar{u}}M(u):=\left\{  v\in
V\ |\ \forall W\in\mathcal{N}_{V}(v),\forall Z\in\mathcal{N}_{U}(\bar
{u}),\text{ }\exists u\in Z,M(u)\cap W\neq\emptyset\right\}  .
\]
Equivalently, $v\in\tau$-$\limsup\limits_{U\ni u\rightarrow\bar{u}}M(u)$ iff
$v$ is the $\tau$-limit of a net $(v_{\alpha})$ such that $v_{\alpha}\in
M(u_{\alpha})$ for some\ $(v_{\alpha})\subset U$ converging to $\bar{u}$. If
the sets $U$ and $V$ are first countable, then we take sequences instead of
nets. We will often omit the reference to $\tau$ and just write $\limsup
_{u\rightarrow\bar{u}}M(u)$ when the topology $\tau$ is understood.

\section{$\varepsilon$-subdifferential calculus for pointwise
suprema\label{sect2}}

In this section, we develop different rules for the $\varepsilon
$-subdifferential mapping of pointwise suprema. The setting here is that of a
dual pair $(X,Y)$ of (real) vector spaces with an associated separating
bilinear form denoted by $\left\langle \cdot,\cdot\right\rangle ,\ $so that
$X$ and $Y$ are endowed with compatible topologies.

In the following we characterize\ the $\varepsilon$-subdifferential mapping of
the conjugate function.

\begin{lem}
\label{thm1}Given a function $f:X\rightarrow\mathbb{R\cup\{+\infty\}}$ such
that $f^{\ast}$ is proper, we have for all $\varepsilon>0$ and $x^{\ast}\in
Y$
\begin{equation}
\partial_{\varepsilon}f^{\ast}(x^{\ast})=\bigcap_{\delta>0}\operatorname*{cl}%
\left(  \sum_{i\in\overline{1,k}}\lambda_{i}(\partial_{\varepsilon_{i}}%
f)^{-1}(x^{\ast})\mid(\lambda_{i})\in\Delta_{k},\text{ }\sum_{i\in
\overline{1,k}}\lambda_{i}\varepsilon_{i}\leq\varepsilon+\delta,\text{
}\varepsilon_{i}\geq0,\text{ }k\geq1\right)  . \label{first}%
\end{equation}

\end{lem}

\begin{dem}
To verify the inclusion \textquotedblleft$\supset$\textquotedblright, we fix
$\delta>0$ and take $x=\sum_{i=1}^{k}\lambda_{i}x_{i}$ with $x_{i}\in
(\partial_{\varepsilon_{i}}f)^{-1}(x^{\ast})$, $(\lambda_{i})\in\Delta_{k}$
and $\sum\lambda_{i}\varepsilon_{i}\leq\varepsilon+\delta$ $(\varepsilon
_{i}\geq0,$ $k\in\mathbb{N}).$ Then $x^{\ast}\in\partial_{\varepsilon_{i}%
}f(x_{i})$ and, so,
\[
\left\langle x^{\ast},u-x_{i}\right\rangle \leq f(u)-f(x_{i})+\varepsilon
_{i},\,\forall u\in X,\text{ }\forall i=1,\cdots,k.
\]
Multiplying this inequality by $\lambda_{i}$ and summing up over $i,$ we
obtain (recall that $f^{\ast\ast}=(f^{\ast})^{\ast}=\overline
{\operatorname*{co}}f$, by (\ref{cob}))
\begin{equation}
\left\langle x^{\ast},u-x\right\rangle \leq f(u)+\sum_{i\in\overline{1,k}%
}\lambda_{i}(-f(x_{i})+\varepsilon_{i})\leq f(u)-f^{\ast\ast}(x)+\varepsilon
+\delta. \label{su}%
\end{equation}
Since $x^{\ast}\in\partial_{\varepsilon_{i}}f(x_{i})$ and $f^{\ast}$ is
proper, we have $-\infty<f^{\ast\ast}(x)=f^{\ast\ast}(\sum_{i=1}^{k}%
\lambda_{i}x_{i})\leq\sum_{i\in\overline{1,k}}\lambda_{i}f^{\ast\ast}%
(x_{i})\leq\sum_{i\in\overline{1,k}}\lambda_{i}f(x_{i})<+\infty.$ Then, by
taking the supremum over $u$ in (\ref{su}),%
\[
f^{\ast}(x^{\ast})+(f^{\ast})^{\ast}(x)=f^{\ast}(x^{\ast})+f^{\ast\ast}%
(x)\leq\left\langle x^{\ast},x\right\rangle +\varepsilon+\delta;
\]
that is, $x\in\partial_{\varepsilon+\delta}f^{\ast}(x^{\ast}).$ Because
$\delta$ was arbitrarily chosen, it follows that $x\in\partial_{\varepsilon
}f^{\ast}(x^{\ast}).$ Consequently, the inclusion \textquotedblleft$\supset
$\textquotedblright\ follows due to the closedness and the convexity of the
set $\partial_{\varepsilon}f^{\ast}(x^{\ast}).$

To establish the\ inclusion \textquotedblleft$\subset$\textquotedblright\ we
take $x\in\partial_{\varepsilon}f^{\ast}(x^{\ast})$ so that, by (\ref{cob}),
$\overline{\operatorname*{co}}f(x)=f^{\ast\ast}(x)\in\mathbb{R}$ and
\begin{equation}
f^{\ast}(x^{\ast})+\overline{\operatorname*{co}}f(x)\leq\left\langle x^{\ast
},x\right\rangle +\varepsilon. \label{asu}%
\end{equation}
Given a $\delta>0$ we choose a $V\in\mathcal{N}_{X}(\theta)$ such that
$\mathrm{\sigma}_{V}(x^{\ast})\leq\delta.$ From\ the definition of
$\overline{\operatorname*{co}}f$ there are elements\ $x_{1},\cdots,x_{k}%
\in\operatorname*{dom}f,$ $(\lambda_{i})\in\Delta_{k},$ and $k\in\mathbb{N},$
such that $x-\sum_{i=1}^{k}\lambda_{i}x_{i}\in V$ and
\begin{equation}
\overline{\operatorname*{co}}f(x)-\delta\leq\sum_{i=1}^{k}\lambda_{i}%
f(x_{i})\leq\overline{\operatorname*{co}}f(x)+\delta\leq-f^{\ast}(x^{\ast
})+\left\langle x^{\ast},x\right\rangle +\varepsilon+\delta. \label{nm}%
\end{equation}
We put\ $\varepsilon_{i}:=f(x_{i})-\overline{\operatorname*{co}}%
f(x)+\left\langle x^{\ast},x-x_{i}\right\rangle +\varepsilon$ $(\in
\mathbb{R}).$ Observe that from\ the definition of $f^{\ast}$ and the relation
(\ref{asu})\ we have\
\[
\varepsilon_{i}=-(\left\langle x^{\ast},x_{i}\right\rangle -f(x_{i}%
))-\overline{\operatorname*{co}}f(x)+\left\langle x^{\ast},x\right\rangle
+\varepsilon\geq-f^{\ast}(x^{\ast})-\overline{\operatorname*{co}%
}f(x)+\left\langle x^{\ast},x\right\rangle +\varepsilon\geq0.
\]
Also, using (\ref{nm}) and the choice of $V$ we obtain (recall (\ref{cob}))%
\begin{align*}
\sum_{i=1}^{k}\lambda_{i}\varepsilon_{i}  &  =\sum_{i=1}^{k}\lambda_{i}%
f(x_{i})-\overline{\operatorname*{co}}f(x)+\left\langle x^{\ast},x-\sum
_{i=1}^{k}\lambda_{i}x_{i}\right\rangle +\varepsilon\\
&  \leq\delta+\mathrm{\sigma}_{V}(x^{\ast})+\varepsilon\leq2\delta
+\varepsilon.
\end{align*}
Thus, since (recall (\ref{asu}))
\[
\varepsilon_{i}+\left\langle x^{\ast},x_{i}\right\rangle =f(x_{i}%
)+(\left\langle x^{\ast},x\right\rangle +\varepsilon-\overline
{\operatorname*{co}}f(x))\geq f(x_{i})+f^{\ast}(x^{\ast}),
\]
it follows that $x_{i}\in(\partial_{\varepsilon_{i}}f)^{-1}(x^{\ast}),$ and
so\ $x\in\sum_{i=1}^{k}\lambda_{i}(\partial_{\varepsilon_{i}}f)^{-1}(x^{\ast
})+V.$ Finally, the arbitrariness of $V$ and $\delta>0$ leads us to\ the
desired inclusion \textquotedblleft$\subset$\textquotedblright.
\end{dem}

\bigskip

We give now a formula for\ the $\varepsilon$-subdifferential of the supremum
function, which extends and improves\ \cite[Theorem 1]{Volle12}.

\begin{thm}
\label{maxmain}Given set $T$ and functions $\Phi_{t}\in\Gamma_{0}(X)$, $t\in
T,$ we put $\Phi:=\sup_{t\in T}\Phi_{t}.$ Then for every $x\in X$ and
$\varepsilon\geq0$ we have
\[
\partial_{\varepsilon}\Phi(x)=\bigcap_{\delta>\varepsilon}\operatorname*{cl}%
\left(  \bigcup\limits_{\substack{(\lambda_{i})\in\Delta_{k},\text{ }t_{i}\in
T,\text{ }\beta_{i}\geq0,\text{ }k\in\mathbb{N}\\\sum\limits_{i\in
\overline{1,k}}\lambda_{i}(\beta_{i}-\Phi_{t_{i}}(x)+\Phi(x))=\delta}%
}\sum_{i\in\overline{1,k}}\lambda_{i}\partial_{\beta_{i}}\Phi_{t_{i}%
}(x)\right)  .
\]

\end{thm}

\begin{dem}
We fix $x\in X$ and $\delta>\varepsilon\geq0.$ If\ $g:=\inf_{t\in T}\Phi
_{t}^{\ast},$ then\ writing (by (\ref{cob}))\
\begin{equation}
g^{\ast}=\sup_{t\in T}(\Phi_{t}^{\ast})^{\ast}=\sup_{t\in T}\overline
{\operatorname*{co}}\Phi_{t}=\sup_{t\in T}\Phi_{t}=\Phi, \label{bi}%
\end{equation}
according to Lemma\ \ref{thm1} it follows that
\begin{align}
\partial_{\varepsilon}\Phi(x)  &  =\partial_{\varepsilon}g^{\ast
}(x)\nonumber\\
&  =\bigcap_{\delta>\varepsilon}\operatorname*{cl}\left(  \sum\limits_{i\in
\overline{1,k}}\lambda_{i}(\partial_{\varepsilon_{i}}g)^{-1}(x)\mid
(\lambda_{i})\in\Delta_{k},\text{ }\sum_{i\in\overline{1,k}}\lambda
_{i}\varepsilon_{i}\leq\delta,\text{ }\varepsilon_{i}\geq0,\text{ }%
k\geq1\right)  . \label{ax}%
\end{align}

To establish the inclusion \textquotedblleft$\subset$\textquotedblright\ of
the current theorem we pick $x_{i}^{\ast}\in(\partial_{\varepsilon_{i}}%
g)^{-1}(x),$ $i\in\overline{1,k},$ where $\varepsilon_{i}\geq0$ and
$k\in\mathbb{N}$ are such that $\sum_{i\in\overline{1,k}}\lambda
_{i}\varepsilon_{i}\leq\delta$ for some $(\lambda_{i})\in\Delta_{k}.$
Then\ $x\in\partial_{\varepsilon_{i}}g(x_{i}^{\ast})$ and
\[
\inf_{t\in T}\Phi_{t}^{\ast}(x_{i}^{\ast})+\Phi(x)=g(x_{i}^{\ast})+g^{\ast
}(x)\leq\left\langle x_{i}^{\ast},x\right\rangle +\varepsilon_{i}.
\]
By choosing $t_{i}\in T$ such that $\inf_{t\in T}\Phi_{t}^{\ast}(x_{i}^{\ast
})\geq\Phi_{t_{i}}^{\ast}(x_{i}^{\ast})-\delta+\varepsilon$ we obtain
\begin{equation}
\Phi_{t_{i}}^{\ast}(x_{i}^{\ast})+\Phi(x)\leq\inf_{t\in T}\Phi_{t}^{\ast
}(x_{i}^{\ast})+\Phi(x)+\delta-\varepsilon\leq\left\langle x_{i}^{\ast
},x\right\rangle +\varepsilon_{i}+\delta-\varepsilon, \label{this}%
\end{equation}
and so $\Phi(x)\leq\left\langle x_{i}^{\ast},x\right\rangle -\Phi_{t_{i}%
}^{\ast}(x_{i}^{\ast})+\varepsilon_{i}+\delta-\varepsilon\leq(\Phi_{t_{i}%
}^{\ast})^{\ast}(x)+\varepsilon_{i}+\delta-\varepsilon=\Phi_{t_{i}%
}(x)+\varepsilon_{i}+\delta-\varepsilon.$ Also, since\ $\Phi_{t_{i}}^{\ast
}(x_{i}^{\ast})+\Phi(x)-\left\langle x_{i}^{\ast},x\right\rangle
\leq\varepsilon_{i}+\delta-\varepsilon,$ we get\
\[
0\leq\Phi_{t_{i}}^{\ast}(x_{i}^{\ast})+\Phi_{t_{i}}(x)-\left\langle
x_{i}^{\ast},x\right\rangle \leq\varepsilon_{i}+\delta-\varepsilon+\Phi
_{t_{i}}(x)-\Phi(x)=:\hat{\beta}_{i};
\]
that is, $x_{i}^{\ast}\in\partial_{\hat{\beta}_{i}}\Phi_{t_{i}}(x)$ and
$\sum\limits_{i\in\overline{1,k}}\lambda_{i}(\hat{\beta}_{i}-\Phi_{t_{i}%
}(x))=\sum\limits_{i\in\overline{1,k}}\lambda_{i}(\varepsilon_{i}%
+\delta-\varepsilon-\Phi(x))\leq2\delta-\varepsilon-\Phi(x).$

Let $\gamma\geq0$ such that $\sum\limits_{i\in\overline{1,k}}\lambda_{i}%
(\hat{\beta}_{i}-\Phi_{t_{i}}(x))+\gamma=2\delta-\varepsilon-\Phi(x)\ $and
denote $\beta_{i}:=\hat{\beta}_{i}+\gamma.$ Then $x_{i}^{\ast}\in
\partial_{\hat{\beta}_{i}}\Phi_{t_{i}}(x)\subset\partial_{\beta_{i}}%
\Phi_{t_{i}}(x)$, $\sum\limits_{i\in\overline{1,k}}\lambda_{i}(\beta_{i}%
-\Phi_{t_{i}}(x))=2\delta-\varepsilon-\Phi(x),\ $and, consequently, the
inclusion \textquotedblleft$\subset$\textquotedblright\ follows thanks
to\ (\ref{ax}).

To prove the inclusion \textquotedblleft$\supset$\textquotedblright\ we
take\ $x^{\ast}=$ $\sum_{i\in\overline{1,k}}\lambda_{i}x_{i}^{\ast},$ where
$x_{i}^{\ast}\in\partial_{\beta_{i}}\Phi_{t_{i}}(x)$ $(i\in\overline{1,k}),$
$\delta>0,$ $k\in\mathbb{N},$ and $(\lambda_{i})\in\Delta_{k},\ $%
with\ $\sum\limits_{i\in\overline{1,k}}\lambda_{i}(\beta_{i}-\Phi_{t_{i}%
}(x))+\Phi(x)=\varepsilon+\delta.$ Then
\[
\Phi_{t_{i}}(x)+\Phi_{t_{i}}^{\ast}(x_{i}^{\ast})\leq\left\langle
x,x_{i}^{\ast}\right\rangle +\beta_{i},
\]
and, by multiplying by $\lambda_{i}$ and summing over $i\in\overline{1,k},$
\[
\sum\limits_{i\in\overline{1,k}}\lambda_{i}\Phi_{t_{i}}^{\ast}(x_{i}^{\ast
})\leq\left\langle x,\sum\limits_{i\in\overline{1,k}}\lambda_{i}x_{i}^{\ast
}\right\rangle +\sum\limits_{i\in\overline{1,k}}\lambda_{i}(\beta_{i}%
-\Phi_{t_{i}}(x))=\left\langle x,x^{\ast}\right\rangle +\varepsilon
+\delta-\Phi(x).
\]
But $\Phi^{\ast}=(g^{\ast})^{\ast}=\overline{\operatorname*{co}}(\inf_{t\in
T}\Phi_{t}^{\ast})$ (as $g^{\ast}=\Phi$ is proper, recall (\ref{bi})), and so
the last inequality yields $\Phi^{\ast}(x^{\ast})\leq\sum\limits_{i\in
\overline{1,k}}\lambda_{i}\overline{\operatorname*{co}}(\inf_{t\in T}\Phi
_{t}^{\ast})(x_{i}^{\ast})\leq\left\langle x,x^{\ast}\right\rangle
+\varepsilon+\delta-\Phi(x);$ that is, $x^{\ast}\in\partial_{\varepsilon
+\delta}\Phi(x).$ Thus, the desired inclusion \textquotedblleft$\supset
$\textquotedblright\ follows along the arbitrariness of $\delta>0.$
\end{dem}

\bigskip

Theorem \ref{maxmain} becomes more explicit when the index set $T$ is
countable and the sequence\ $(\Phi_{n})_{n}$ is non-decreasing.

\begin{cor}
\label{increasing0}With the assumptions of Theorem \emph{\ref{maxmain}}\ we
take $T=\mathbb{N}.$ If\ the sequence $(\Phi_{n})_{n}$ is non-decreasing, then
for every $x\in X$ and $\varepsilon\geq0$ we have
\[
\partial_{\varepsilon}\Phi(x)=\bigcap_{\delta>\varepsilon}\limsup
_{n\rightarrow+\infty}\partial_{\delta}\Phi_{n}(x)=\limsup_{n\rightarrow
+\infty,\text{ }\delta\downarrow\varepsilon}\partial_{\delta}\Phi_{n}(x).
\]

\end{cor}

\begin{dem}
Take $\xi\in\partial_{\varepsilon}\Phi(x)$ and fix $\delta>0,$ $V\in
\mathcal{N}_{Y}(\theta).$ According to Theorem \ref{maxmain}, we have that
$\xi\in\sum\limits_{i=\overline{1,k}}\lambda_{i}\xi_{i}+V$ for some $\xi
_{i}\in\partial_{\beta_{i}}\Phi_{n_{i}}(x)$ and $(\lambda_{i})\in\Delta_{k}$
$(k\in\mathbb{N)},$ where $\beta_{i}\geq0$ and $n_{i}\in\mathbb{N}%
\ (i\in\overline{1,k})$ are such that $\sum\limits_{i=\overline{1,k}}%
\lambda_{i}(\beta_{i}-\Phi_{n_{i}}(x)+\Phi(x))=\varepsilon+\frac{\delta}{2}.$
Set\ $m_{0}:=\max_{i\in\overline{1,k}}n_{i}\geq1.$ Then, from one hand,
by\ the current assumption on the sequence $(\Phi_{n}),$ we obtain
\begin{equation}
\Phi_{m_{0}}(x)\geq\sum\limits_{i=\overline{1,k}}\lambda_{i}\Phi_{n_{i}%
}(x)\geq\Phi(x)+\sum\limits_{i=\overline{1,k}}\lambda_{i}\beta_{i}%
-\varepsilon-\frac{\delta}{2}\geq\Phi(x)-\varepsilon-\delta. \label{ie}%
\end{equation}
On the other hand, by writing the relation\ $\xi_{i}\in\partial_{\beta_{i}%
}\Phi_{n_{i}}(x)$ into an inequality form\ and, next, summing up over
$i\in\overline{1,k},$ we obtain, for all $u\in x+\frac{\delta}{2}V^{\circ}$
(hence, $\mathrm{\sigma}_{V}(u-x)\leq\frac{\delta}{2}),$
\begin{align*}
\left\langle \xi,u-x\right\rangle  &  \leq\sum\limits_{i=\overline{1,k}%
}\lambda_{i}\left\langle \xi_{i},u-x\right\rangle +\frac{\delta}{2}\\
&  \leq\sum\limits_{i=\overline{1,k}}\lambda_{i}(\Phi_{n_{i}}(u)-\Phi_{n_{i}%
}(x)+\beta_{i})+\frac{\delta}{2}\\
&  \leq\sum\limits_{i=\overline{1,k}}\lambda_{i}(\Phi_{n_{i}}(u)-\Phi
(x))+\varepsilon+\delta\\
&  \leq\Phi_{m_{0}}(u)-\Phi(x)+\varepsilon+\delta.
\end{align*}
Thus,
\[
\left\langle \xi,u-x\right\rangle \leq\Phi_{n}(y)-\Phi_{n}(x)+\varepsilon
+\delta\text{ \ \ for all }n\geq m_{0},
\]
and, so, taking into account (\ref{ie}), by the sum rule of $\varepsilon
$-subdifferentials (e.g., \cite{JBHUPH}) we get
\[
\xi\in%
{\textstyle\bigcap\limits_{n\geq m_{0}}}
\partial_{\varepsilon+\delta}(\Phi_{n}+\mathrm{I}_{x+\frac{\delta}{2}V^{\circ
}})(x)\subset%
{\textstyle\bigcap\limits_{n\geq m_{0},\text{ }\eta>0}}
\partial_{\varepsilon+\delta+\eta}\Phi_{n}(x)+\frac{3(\varepsilon+\delta
+\eta)}{\delta}V.
\]
Hence, as $V$ and $\delta$ were arbitrarily chosen, we deduce that
\[
\xi\in\limsup_{n\rightarrow+\infty,\text{ }\delta\downarrow\varepsilon
}\partial_{\varepsilon}\Phi_{n}(x)=\bigcap_{\delta>0}\limsup_{n\rightarrow
+\infty}\partial_{\varepsilon+\delta}\Phi_{n}(x).
\]
This finishes the proof since the opposite inclusion \textquotedblleft%
$\supset$\textquotedblright\ holds straightforwardly.
\end{dem}

\bigskip

We recover the subdifferenial rule for\ the case of finitely many convex
functions; see, e.g., \cite[Corollary 2.8.11]{zalinesku2000convex}.

\begin{cor}
\label{finitecase}With the assumptions of Theorem \emph{\ref{maxmain}}\ we
take $T=\{1,\cdots,n\}$. Then for every $x\in X$%
\begin{align*}
\partial_{\varepsilon}\Phi(x)=\bigcup\{\partial_{\eta}\left(  \sum
\nolimits_{i\in\overline{1,n}}\lambda_{i}\Phi_{i}(x)\right)  (x)\mid\text{ }
&  (\lambda_{i})\in\Delta_{n},\text{ }\\
&  \eta\in\left[  0,\varepsilon\right]  ,\text{ }\sum\nolimits_{i=1}%
^{n}\lambda_{i}\Phi_{i}(x)\geq\Phi(x)+\eta-\varepsilon\}.
\end{align*}

\end{cor}

\begin{dem}
We may assume that $x\in\operatorname*{dom}\Phi=\cap_{i\in\overline{1,n}%
}\operatorname*{dom}\Phi_{i}$ and $\partial_{\varepsilon}\Phi(x)\neq
\emptyset.$ The inclusion \textquotedblleft$\supset$\textquotedblright\ is
straightforward. For the other inclusion we have, by Theorem \ref{maxmain},
\begin{align*}
\partial_{\varepsilon}\Phi(x)  &  =\bigcap_{\delta>\varepsilon}%
\operatorname*{cl}\{\sum\nolimits_{i=1}^{n}\lambda_{i}\partial_{\beta_{i}}%
\Phi_{i}(x)\mid(\lambda_{i})\in\Delta_{n},\ \beta_{i}\geq0,\text{ }\\
&  \text{~\quad\quad\quad\quad\quad\quad\quad\quad\quad\quad\quad\quad
\quad\quad\quad\quad\quad\quad}\sum\nolimits_{i=1}^{n}\lambda_{i}(\beta
_{i}-\Phi_{i}(x))+\Phi(x)=\delta\}\\
&  =\bigcap_{\delta>\varepsilon}\operatorname*{cl}\{\sum\nolimits_{i=1}%
^{n}\lambda_{i}\partial_{\beta_{i}}\Phi_{i}(x)\mid(\lambda_{i})\in\Delta
_{n},\ \lambda_{i},\text{ }\beta_{i}>0,\text{ }\\
&  \text{~\quad\quad\quad\quad\quad\quad\quad\quad\quad\quad\quad\quad
\quad\quad\quad\quad\quad\quad}\sum\nolimits_{i=1}^{n}\lambda_{i}(\beta
_{i}-\Phi_{i}(x))+\Phi(x)=\delta\},
\end{align*}
and so, setting $\eta=\sum_{i=1}^{n}\lambda_{i}\beta_{i},$
\begin{align}
\partial_{\varepsilon}\Phi(x)  &  \subset\bigcap_{\delta>\varepsilon
}\operatorname*{cl}\{\partial_{\eta}\left(  \sum\nolimits_{i=1}^{n}\lambda
_{i}\Phi_{i}\right)  (x)\mid(\lambda_{i})\in\Delta_{n},\ \eta\in\left[
0,\varepsilon\right]  ,\text{ }\nonumber\\
&  \text{~\quad\quad\quad\quad\quad\quad\quad\quad\quad\quad\quad\quad
\quad\quad\quad\quad\quad\quad}\sum\nolimits_{i=1}^{n}\lambda_{i}\Phi
_{i}(x)\geq\Phi(x)+\eta-\delta\}. \label{oi}%
\end{align}
Then the desired inclusion\ follows due\ to the compactness of the set
$\Delta_{n}.$
\end{dem}

\section{$\varepsilon$-subdifferential approach\label{sect3}}

In this section, we give the desired characterization of $\delta$-normal sets
to arbitrary intersections of sublevel sets. As in the previous section, the
framework here is that of\ a dual pair $(X,Y)$ of (real) lc spaces $X$ and
$Y,$ endowed with compatible topologies with respect to a given dual pairing
$\left\langle \cdot,\cdot\right\rangle $. We consider\ a family\ of proper,
lsc, and convex functions
\[
\Phi_{t}:X\rightarrow\mathbb{R\cup\{+\infty\}},\text{ \ }t\in T,
\]
where $T$ is an arbitrary index set, together with the associated supremum
function
\[
\Phi:=\sup_{t\in T}\Phi_{t}.
\]
The following theorem provides the main result of this section. Its proof is
based on a series of lemmas that we postpone to the Appendix.

\begin{thm}
\label{tmain}For every $\bar{x}\in\operatorname*{dom}\Phi,$ $\delta\geq0$ and
$\lambda\in\left(  -\infty,+\infty\right]  $ we have
\begin{equation}
\mathrm{N}_{([\Phi\leq\lambda]\cap\operatorname*{dom}\Phi)\cup(\bar{x}%
+[\Phi\leq\Phi(\bar{x})]_{\infty})}^{\delta}(\bar{x})=\limsup_{\substack{\sum
\nolimits_{i\in\overline{1,k}}\mu_{i}(\nu-\Phi_{t_{i}}(\bar{x})+\varepsilon
_{i})\rightarrow\delta,\text{ }\nu\uparrow\lambda\\\mu_{i}>0,\text{
}\varepsilon_{i}\geq0,\text{ }t_{i}\in T,\text{ }k\in\mathbb{N}}%
}\sum\nolimits_{i\in\overline{1,k}}\mu_{i}\partial_{\varepsilon_{i}}%
\Phi_{t_{i}}(\bar{x}). \label{generalformula}%
\end{equation}

\end{thm}

\begin{dem}
Let us start with the proof of the inclusion \textquotedblleft$\supset
$\textquotedblright. We take $x^{\ast}=\lim_{j}\sum_{i=1}^{k_{j}}\mu
_{i,j}x_{i,j}^{\ast}$ for $k_{j}\in\mathbb{N},$ $\mu_{i,j}>0,$ $\varepsilon
_{i,j}\geq0,$ $t_{i,j}\in T$ and $x_{i,j}^{\ast}\in\partial_{\varepsilon
_{i,j}}\Phi_{t_{i,j}}(\bar{x})$ $(i=\overline{1,k_{j}})$ such that%
\[
\sum_{i=0}^{k_{j}}\mu_{i,j}(\nu_{j}-\Phi_{t_{i,j}}(\bar{x})+\varepsilon
_{i,j})\rightarrow\delta,
\]
where\ $\nu_{j}\uparrow\lambda$ (observe that one can take $\nu_{j}=\lambda$
when $\lambda$ is finite).

We fix an element\ $u\in\lbrack\Phi\leq\lambda]\cap\operatorname*{dom}\Phi;$
hence, we may suppose that $\Phi(u)\leq\nu_{j}$ for all $j.$ Then from\ the
definition of $x_{i,j}^{\ast}$ we get
\begin{equation}
\left\langle x_{i,j}^{\ast},u-\bar{x}\right\rangle \leq\Phi_{t_{i,j}}%
(u)-\Phi_{t_{i,j}}(\bar{x})+\varepsilon_{i,j}\leq\nu_{j}-\Phi_{t_{i,j}}%
(\bar{x})+\varepsilon_{i,j}. \label{ineqsubdarbi}%
\end{equation}
Multiplying this last inequality by $\mu_{i,j}$ and summing up\ over
$i=\overline{1,k_{j}}$ give us
\[
\left\langle \sum_{i=1}^{k_{j}}\mu_{i,j}x_{i,j}^{\ast},u-\bar{x}\right\rangle
\leq\sum_{i=0}^{k_{j}}\mu_{i,j}(\nu_{j}-\Phi_{t_{i,j}}(\bar{x})+\varepsilon
_{i,j}),
\]
which at the limit yields
\begin{equation}
\left\langle x^{\ast},u-\bar{x}\right\rangle \leq\delta. \label{ry}%
\end{equation}
Now we take $u\in\lbrack\Phi\leq\Phi(\bar{x})]_{\infty}.$ Then for all
$\beta>0$ we have $\bar{x}+\beta u\in\bar{x}+[\Phi\leq\Phi(\bar{x})]_{\infty}$
and, so, by (\ref{ineqsubdarbi}),
\[
\left\langle \sum_{i=1}^{k_{j}}\mu_{i,j}x_{i,j}^{\ast},\beta u\right\rangle
=\left\langle \sum_{i=1}^{k_{j}}\mu_{i,j}x_{i,j}^{\ast},\bar{x}+\beta
u-\bar{x}\right\rangle \leq\sum_{i=0}^{k_{j}}\mu_{i,j}(\Phi(\bar{x}%
)-\Phi_{t_{i,j}}(\bar{x})+\varepsilon_{i,j}),
\]
implying that $\left\langle \sum_{i=1}^{k_{j}}\mu_{i,j}x_{i,j}^{\ast
},u\right\rangle \leq0.$ Hence, after taking the the limit on $j$ we
get\ $\left\langle x^{\ast},u\right\rangle \leq0\leq\delta.$ Combining this
with (\ref{ry}) yields\ $x^{\ast}\in\mathrm{N}_{([\Phi\leq\lambda
]\cap\operatorname*{dom}\Phi)\cup(\bar{x}+[\Phi\leq\Phi(\bar{x})]_{\infty}%
)}^{\delta}(\bar{x}).$ This proves the inclusion \textquotedblleft$\supset
$\textquotedblright.

To prove the inclusion \textquotedblleft$\subset$\textquotedblright\ we pick
an element\ $\xi\in\mathrm{N}_{([\Phi\leq\lambda]\cap\operatorname*{dom}%
\Phi)\cup(\bar{x}+[\Phi\leq\Phi(\bar{x})]_{\infty})}^{\delta}(\bar{x}).$ We
proceed by investigating all the possible values of $\lambda$ and $\Phi
(\bar{x}):$

$(1)$ $\Phi(\bar{x})\leq\lambda<+\infty.$ In this case, the right-hand side in
(\ref{generalformula}) reduces to $\mathrm{N}_{[\Phi\leq\lambda]}^{\delta
}(\bar{x}),\ $and it follows from Lemma \ref{lem3} that $\xi\in\limsup
_{\mu(\lambda-\Phi(\bar{x})+\varepsilon)\rightarrow\delta,\text{ }\mu>0}%
\mu\partial_{\varepsilon}\Phi(\bar{x}).$ Then, for every given\ $\eta>0$
(small enough) and\ $\theta$-neighborhood $V\subset Y$ $(V\in\mathcal{N}%
_{Y}(\theta)),$ there exist $\mu>0,$ $\varepsilon\geq0$ and $x^{\ast}%
\in\partial_{\varepsilon}\Phi(\bar{x})$ such that
\begin{equation}
\xi\in\mu x^{\ast}+\frac{1}{2}V\text{, }\left\vert \mu(\lambda-\Phi(\bar
{x})+\varepsilon)-\delta\right\vert <\frac{\eta}{2}. \label{gf}%
\end{equation}
Thus,\ by Theorem \ref{maxmain} there are $k\in\mathbb{N}$, $t_{i}\in T,$
$\varepsilon_{i}\geq0$ and $x_{i}^{\ast}\in\partial_{\varepsilon_{i}}%
\Phi_{t_{i}}(\bar{x})$ $(i\in\overline{1,k})$, together with $(\alpha_{i}%
)\in\Delta_{k}$ such that
\[
\varepsilon-\frac{\eta}{2\mu}\leq\sum_{i=1}^{k}\alpha_{i}(\varepsilon_{i}%
-\Phi_{t_{i}}(\bar{x}))+\Phi(\bar{x})\leq\varepsilon+\frac{\eta}{2\mu},\text{
}x^{\ast}\in\sum_{i=1}^{k}\alpha_{i}x_{i}^{\ast}+\frac{1}{2\mu}V.
\]
Therefore, for $\mu_{i}:=\mu\alpha_{i}$ $(>0)$ we obtain that
\[
\sum_{i=1}^{k}\mu_{i}(\lambda-\Phi_{t_{i}}(\bar{x})+\varepsilon_{i})=\mu
\sum_{i=1}^{k}\alpha_{i}(\varepsilon_{i}-\Phi_{t_{i}}(\bar{x})+\lambda)\leq
\mu(\lambda-\Phi(\bar{x})+\varepsilon)+\frac{\eta}{2}<\delta+\eta,
\]%
\[
\sum_{i=1}^{k}\mu_{i}(\lambda-\Phi_{t_{i}}(\bar{x})+\varepsilon_{i})=\mu
\sum_{i=1}^{k}\alpha_{i}(\varepsilon_{i}-\Phi_{t_{i}}(\bar{x})+\lambda)\geq
\mu(\lambda-\Phi(\bar{x})+\varepsilon)-\frac{\eta}{2}>\delta-\eta,
\]
together with $\xi\in\sum_{i=1}^{k}\mu_{i}x_{i}^{\ast}+V$. This yields the
inclusion \textquotedblleft$\subset$\textquotedblright\ in
(\ref{generalformula}).

$(2)$ $\lambda=+\infty.$ In this case, the right-hand side in
(\ref{generalformula}) reduces to $\mathrm{N}_{\operatorname*{dom}\Phi
}^{\delta}(\bar{x})$ and so, from Lemma \ref{lem2b}, for every given\ $\eta>0$
(small enough) and\ $\theta$-neighborhood $V\subset Y\ $there exist $\mu
\in(0,\frac{\eta^{2}}{2})$ and $\varepsilon\geq0$ together with $x^{\ast}%
\in\partial_{\varepsilon}\Phi(\bar{x})$ such that $\mu\varepsilon\in
(\delta-\frac{\eta}{2},\delta+\frac{\eta}{2})$ and\
\[
\xi\in\mu x^{\ast}+V.
\]
We take $\nu=\eta^{-1}+\Phi(\bar{x})$ $(>0),$ so that $\mu\nu=\frac{\eta}{2}$
and $\left\vert \mu(\nu-\Phi(\bar{x})+\varepsilon)-\delta\right\vert
\leq\left\vert \mu\varepsilon-\delta\right\vert +\left\vert \mu(\nu-\Phi
(\bar{x}))\right\vert <\frac{\eta}{2}+\frac{\eta}{2}=\eta.$ In other words, we
also have (\ref{gf}) in the current case, and the proof follows by arguing as
in point (1) above.

$(3)$\ $\Phi(\bar{x})>\lambda:$ In this case we appeal to\ Lemma \ref{lem4},
which ensures that $\xi\in\overline{%
{\textstyle\bigcup\limits_{\mu>0}}
\partial_{\delta+\mu(\Phi(\bar{x})-\lambda)}(\mu\Phi)(\bar{x})}=\overline{%
{\textstyle\bigcup\limits_{\mu>0}}
\mu\partial_{\frac{\delta}{\mu}+\Phi(\bar{x})-\lambda}\Phi(\bar{x})}.$ Hence,
for every $\theta$-neighborhood $V\subset Y,$ there exist $\mu>0\ $and
$x^{\ast}\in\partial_{\frac{\delta}{\mu}+\Phi(\bar{x})-\lambda}\Phi(\bar{x}%
)$\ such that $\xi\in\mu x^{\ast}+\frac{1}{2}V.$ Hence, (\ref{gf}) follows by
taking $\varepsilon=\frac{\delta}{\mu}+\Phi(\bar{x})-\lambda,$ and we proceed
as in point (1) above.

The proof of the theorem is complete.
\end{dem}

\bigskip

Let us say some words to explain the elements involved in formula
(\ref{generalformula}); namely, the appealing to the set $[\Phi\leq\Phi
(\bar{x})]_{\infty},$ and the consideration of the value $\lambda=+\infty:$

\begin{rem}
$(i)$\emph{ It is clear that when }$\Phi(\bar{x})\leq\lambda<+\infty,$\emph{
the set }$[\Phi\leq\lambda]\cup(\bar{x}+[\Phi\leq\Phi(\bar{x})]_{\infty}%
)$\emph{ reduces to the sublevel set }$[\Phi\leq\lambda],$\emph{ and so
(\ref{generalformula}) gives the required explicit characterization for the
normal cone to }$[\Phi\leq\lambda].$\emph{ The resulting formula in this case
can be compared, when all the }$\Phi_{t}$\emph{'s are all equal, to the one
given in \cite{hiriart1995subdifferential} (see, also, \cite{kutateladze}).
Original characterizations for the normal cone to sublevel sets remounts to
\cite{rockafellar2015convex}. }

$(ii)$\emph{ Even with the lack of the nonemptiness of the set }$[\Phi
\leq\lambda],$\emph{ formula (\ref{generalformula}) is still meaningful, since
the vector }$\bar{x}$ \emph{always belongs to }$\bar{x}+[\Phi\leq\Phi(\bar
{x})]_{\infty}.$ \emph{The other interesting situation covered by
(\ref{generalformula}) is that when }$[\Phi\leq\lambda]$\emph{ is non-empty,
but }$\bar{x}\not \in \lbrack\Phi\leq\lambda].$\emph{ In this case, the
presence of the term }$\bar{x}+[\Phi\leq\Phi(\bar{x})]_{\infty}$\emph{ becomes
necessary, since, for otherwise, the normal cone to }$[\Phi\leq\lambda]$\emph{
at }$\bar{x}$\emph{ can not be defined appropriately.}

$(iii)$\emph{ If }$\lambda=+\infty,$\emph{ then }$[\Phi\leq\lambda
]\cap\operatorname*{dom}\Phi=\operatorname*{dom}\Phi,$\emph{ and the left-hand
side in (\ref{generalformula}) reduces to the normal cone to the domain of
}$\Phi$\emph{. In this case, the relation }$\lambda_{i}\uparrow\lambda$\emph{
is used to force the term }$\sum\nolimits_{i\in\overline{1,k}}\mu_{i}$
\emph{to go to }$+\infty.$\emph{ }

$(iv)$\emph{ Due to the relation}\
\[
\lbrack\Phi\leq\lambda]\cup(\bar{x}+[\Phi\leq\Phi(\bar{x})]_{\infty}%
)=\bigcap\nolimits_{t\in T}([\Phi_{t}\leq\lambda]\cup(\bar{x}+[\Phi_{t}%
\leq\Phi_{t}(\bar{x})]_{\infty})),
\]
\emph{formula\ (\ref{generalformula}) is indeed a characterization of the
normal cone to the arbitrary intersection }$\cap_{t\in T}[\Phi_{t}\leq
\lambda]\cup(\bar{x}+[\Phi_{t}\leq\Phi_{t}(\bar{x})]_{\infty}).$
\end{rem}

We are going now to specify Theorem \ref{tmain} to certain special cases,
which lead to simpler characterizations of the $\delta$-normal set.

Firstly, write formula (\ref{generalformula}) in its most frequent\ form,
corresponding to $\lambda=\Phi(\bar{x}).$

\begin{cor}
\label{cor0}For every $\bar{x}\in\operatorname*{dom}\Phi$ and $\delta\geq0$ we
have that
\[
\mathrm{N}_{[\Phi\leq\Phi(\bar{x})]}^{\delta}(\bar{x})=\limsup_{\substack{\sum
\nolimits_{i\in\overline{1,k}}\mu_{i}(\Phi(\bar{x})-\Phi_{t_{i}}(\bar
{x})+\varepsilon_{i})\rightarrow\delta\\\mu_{i}\geq0,\text{ }\lambda
_{i}\uparrow\lambda,\text{ }t_{i}\in T,\text{ }k\in\mathbb{N}}}\sum
\nolimits_{i\in\overline{1,k}}\mu_{i}\partial_{\varepsilon_{i}}\Phi_{t_{i}%
}(\bar{x}),
\]
and, particularly, when $\Phi\equiv\Phi_{t}$ for all $t\in T,$%
\[
\mathrm{N}_{[\Phi\leq\Phi(\bar{x})]}^{\delta}(\bar{x})=\limsup_{\mu
\varepsilon\rightarrow\delta,\text{ }\mu\geq0}\mu\partial_{\varepsilon}%
\Phi(\bar{x}).
\]

\end{cor}

\begin{dem}
This follows easily from Theorem \ref{tmain} due to the relation $[\Phi
\leq\Phi(\bar{x})]=[\Phi\leq\Phi(\bar{x})]\cup(\bar{x}+[\Phi\leq\Phi(\bar
{x})]_{\infty}).$
\end{dem}

\bigskip

Formula (\ref{generalformula}) takes an algebraic form in the following
corollary, giving rise to a known formula (\cite{kutateladze}).

\begin{cor}
\label{cor1}For every $\bar{x}\in\operatorname*{dom}\Phi$ and $\delta\geq0$ we
have that
\[
\mathrm{N}_{[\Phi\leq\Phi(\bar{x})]}^{\delta}(\bar{x})=\operatorname*{cl}%
\left(  \bigcup_{\mu>0}\partial_{\delta}(\mu\Phi)(\bar{x})\right)  \text{ for
all }\delta>0,
\]
and, consequently,%

\[
\mathrm{N}_{[\Phi\leq\Phi(\bar{x})]}(\bar{x})=\bigcap\nolimits_{\delta
>0}\operatorname*{cl}\left(  \bigcup_{\mu>0}\partial_{\delta}(\mu\Phi)(\bar
{x})\right)  .
\]

\end{cor}

\begin{dem}
Since $[\Phi\leq\Phi(\bar{x})]=[\Phi\leq\lambda]\cup(\bar{x}+[\Phi\leq
\Phi(\bar{x})]_{\infty}),$ by Theorem \ref{tmain} we obtain for all
$\delta>0$
\begin{align*}
\mathrm{N}_{[\Phi\leq\Phi(\bar{x})]}^{\delta}(\bar{x})  &  =\limsup
_{\mu\varepsilon\rightarrow\delta,\text{ }\mu>0}\mu\partial_{\varepsilon}%
\Phi(\bar{x})\\
&  =\limsup_{\mu\varepsilon\rightarrow\delta,\text{ }\mu>0}\partial
_{\mu\varepsilon}(\mu\Phi)(\bar{x})\\
&  \subset\limsup_{\nu\rightarrow\delta}\bigcup_{\mu>0}\partial_{\nu}%
(\frac{\nu}{\delta}\mu\Phi)(\bar{x})\text{ \ (taking }\nu=\mu\varepsilon
\text{)}\\
&  =\limsup_{\nu\rightarrow\delta}\frac{\nu}{\delta}\left(  \bigcup_{\mu
>0}\partial_{\delta}(\mu\Phi)(\bar{x})\right) \\
&  =\operatorname*{cl}\left(  \bigcup_{\mu>0}\partial_{\delta}(\mu\Phi
)(\bar{x})\right) \\
&  \subset\limsup_{\mu\varepsilon\rightarrow\delta,\text{ }\mu>0}\partial
_{\mu\varepsilon}(\mu\Phi)(\bar{x})\text{ (taking }\varepsilon=\frac{\delta
}{\mu}\text{).}%
\end{align*}

\end{dem}

\bigskip

When $[\Phi\leq\lambda]\neq\emptyset$ but $\Phi(\bar{x})>\lambda$ Theorem
\ref{tmain} simplifies to:

\begin{cor}
\label{galb}Given $\bar{x}\in\operatorname*{dom}\Phi\ $and $\lambda
\in\mathbb{R}$ we assume that $[\Phi\leq\lambda]\neq\emptyset$ and $\Phi
(\bar{x})>\lambda.$ Then we have%
\[
\mathrm{N}_{[\Phi\leq\lambda]\cup(\bar{x}+[\Phi\leq\Phi(\bar{x})]_{\infty}%
)}(\bar{x})=\overline{\mathbb{R}_{+}\partial_{\Phi(\bar{x})-\lambda}\Phi
(\bar{x})}.
\]

\end{cor}

\begin{dem}
The inclusion \textquotedblleft$\subset$\textquotedblright\ is immediate from
Lemma \ref{lem2} (see the Appendix), while the converse inclusion follows from
Theorem \ref{tmain}.
\end{dem}

\bigskip

In the following corollary we consider the case in which the sublevel set
$\left[  \Phi\leq\lambda\right]  $ is empty.

\begin{cor}
\label{biz}If $\lambda\in\mathbb{R}$ is such that $[\Phi\leq\lambda
]=\emptyset,$ then for every $\bar{x}\in\operatorname*{dom}\Phi\ $we have
\begin{align*}
\mathrm{N}_{\bar{x}+[\Phi\leq\Phi(\bar{x})]_{\infty}}(\bar{x})  &
=\limsup_{\substack{\mu(\lambda-\Phi(\bar{x})+\varepsilon)\rightarrow
0\\\mu\geq0}}\mu\partial_{\varepsilon}\Phi(\bar{x})\\
&  =\overline{\mathbb{R}_{+}\partial_{\varepsilon}\Phi(\bar{x})}\text{ (for
every }\varepsilon\geq\Phi(\bar{x})-\lambda).
\end{align*}

\end{cor}

\begin{dem}
The first\ equality is a direct consequence\ Theorem \ref{tmain}, while the
inclusion \textquotedblleft$\supset$\textquotedblright\ in the second equality
holds easily. To show the converse inclusion, we fix $\varepsilon\geq\Phi
(\bar{x})-\lambda$\ and take $v\in(\partial_{\varepsilon}\Phi(\bar{x}))^{-}$.
Since $[\Phi\leq\Phi(\bar{x})-\varepsilon]\subset\lbrack\Phi\leq
\lambda]=\emptyset$ we deduce that
\begin{equation}
0\leq\inf_{t>0}\frac{\Phi(\bar{x}+tv)-\Phi(\bar{x})+\varepsilon}{t}%
=\Phi_{\varepsilon}^{\prime}(\bar{x};v)=\sup_{x^{\ast}\in\partial
_{\varepsilon}\Phi(\bar{x})}\left\langle v,x^{\ast}\right\rangle \leq0.
\label{infimum1}%
\end{equation}
Let $t_{n}$ be a minimizing sequence for this last infimum. If $\bar{t}%
\in\mathbb{R}$ is an accumulation point of $t_{n},$ then the last relation
above gives
\[
\Phi(\bar{x}+\bar{t}v)-\Phi(\bar{x})+\varepsilon\leq\liminf_{n}\Phi(\bar
{x}+t_{n}v)-\Phi(\bar{x})+\varepsilon\leq0,
\]
so that $\bar{x}+\bar{t}v\in\lbrack\Phi\leq\lambda]=\emptyset$, a
contradiction. Thus we may assume that $t_{n}\rightarrow+\infty$, so that, by
the convexity of $\Phi,$
\begin{align*}
\sup_{t>0}\frac{\Phi(\bar{x}+tv)-\Phi(\bar{x})}{t}  &  =\lim_{n}\frac
{\Phi(\bar{x}+t_{n}v)-\Phi(\bar{x})}{t_{n}}\\
&  =\lim_{n}\frac{\Phi(\bar{x}+t_{n}v)-\Phi(\bar{x})+\varepsilon}{t_{n}}\\
&  =\inf_{t>0}\frac{\Phi(\bar{x}+tv)-\Phi(\bar{x})+\varepsilon}{t}=0.
\end{align*}
This shows that $\bar{x}+tv\in\lbrack\Phi\leq\Phi(\bar{x})]$ for all $t>0$,
and then $v\in(\mathrm{N}_{\bar{x}+[\Phi\leq\Phi(\bar{x})]_{\infty}}(\bar
{x}))^{-}$. In other words, $(\partial_{\varepsilon}\Phi(\bar{x}))^{-}%
\subset(\mathrm{N}_{\bar{x}+[\Phi\leq\Phi(\bar{x})]_{\infty}}(\bar{x}))^{-}$
and, so, using the bipolar Theorem, $\mathrm{N}_{\bar{x}+[\Phi\leq\Phi(\bar
{x})]_{\infty}}(\bar{x})\subset\left[  (\partial_{\varepsilon}\Phi(\bar
{x}))^{-}\right]  ^{-}\subset\overline{\mathbb{R}_{+}\partial_{\varepsilon
}\Phi(\bar{x})}.$
\end{dem}

The following result\ puts in clear\ the different sets composing the
right-hand side of (\ref{generalformula}).

\begin{cor}
\label{corolarioimportante}For every $\bar{x}\in\operatorname*{dom}\Phi$ and
$\delta\geq0$ we have
\[
\mathrm{N}_{[\Phi\leq\Phi(\overline{x})]}^{\delta}(\overline{x})=\bigcup
_{\mu\geq0}\partial_{\delta}(\mu\Phi)(\overline{x})\cup\limsup_{\mu
\varepsilon\rightarrow\delta,\text{ }\mu\rightarrow+\infty}\mu\partial
_{\varepsilon}\Phi(\overline{x}),
\]
and, in particular,
\[
\mathrm{N}_{[\Phi\leq\Phi(\overline{x})]}(\overline{x})=\mathbb{R}_{+}%
\partial\Phi(\overline{x})\cup\mathrm{N}_{\operatorname*{dom}\Phi}%
(\overline{x})\cup\limsup_{\mu\varepsilon\rightarrow0,\text{ }\mu
\rightarrow+\infty}\mu\partial_{\varepsilon}\Phi(\overline{x}).
\]
Moreover, if Slater's condition holds at $\Phi(\overline{x}),$ then
\[
\mathrm{N}_{[\Phi\leq\Phi(\overline{x})]}^{\delta}(\overline{x})=\bigcup
_{\mu\geq0}\partial_{\delta}(\mu\Phi)(\overline{x}).
\]

\end{cor}

\begin{dem}
To prove the first statement of the corollary we take $x^{\ast}\in
\mathrm{N}_{[\Phi\leq\Phi(\overline{x})]}^{\delta}(\overline{x})$ and,
according to Theorem \ref{tmain}, we let $\mu_{i},$ $\varepsilon_{i}\geq0$ and
$x_{i}^{\ast}\in\partial_{\varepsilon_{i}}\Phi(\overline{x})$ such that
$\mu_{i}\varepsilon_{i}\rightarrow\delta$ and $\mu_{i}x_{i}^{\ast
}\rightharpoonup x^{\ast}.$ We may suppose without loss of generality that
$\mu_{i}$ converges in $\mathbb{R}_{+}\cup\left\{  \infty\right\}  .$ If
$\mu_{i}\rightarrow\mu\in\mathbb{R},$ by taking the limit on $i$ in the
inequality
\[
\left\langle \mu_{i}x_{i}^{\ast},y-\overline{x}\right\rangle \leq(\mu_{i}%
\Phi)(y)-(\mu_{i}\Phi)(\overline{x})+\mu_{i}\varepsilon_{i}\text{ \ for all
}y\in X,\text{ }%
\]
we get $\left\langle x^{\ast},y-\overline{x}\right\rangle \leq(\mu
\Phi)(y)-(\mu\Phi)(\overline{x})+\delta.$ So, $x^{\ast}\in\partial_{\delta
}(\mu\Phi)(\overline{x})$ and the inclusion \textquotedblleft$\subset
$\textquotedblright\ in the first equality holds. The converse inclusion
follows by Theorem \ref{tmain}. For the second formula it suffices to\ observe
that, when $\delta=0$,\
\[
\bigcup_{\mu\geq0}\partial(\mu\Phi)(\overline{x})=\mathbb{R}_{+}\partial
\Phi(\overline{x})\cup\mathrm{N}_{\operatorname*{dom}\Phi}(\bar{x}).
\]

To prove the last statement of the corollary we only need to verify that the
set $\limsup_{\mu\varepsilon\rightarrow\delta,\text{ }\mu\rightarrow\infty}%
\mu\partial_{\varepsilon}\Phi(\overline{x})$ is empty whenever Slater's
condition holds at $\Phi(\overline{x}).$ Otherwise, if this were not the case,
then there would exist $\mu_{i},$ $\varepsilon_{i}\geq0$ such that $\mu
_{i}\varepsilon_{i}\rightarrow\delta$, $\mu_{i}\rightarrow+\infty,$ and
$\mu_{i}x_{i}^{\ast}\rightharpoonup x^{\ast}.$ Then $\varepsilon
_{i}\rightarrow0$ and $x_{i}^{\ast}\rightharpoonup0,$ so that $\theta
\in\limsup_{\varepsilon\rightarrow0}\partial_{\varepsilon}\Phi(\overline
{x})=\partial\Phi(\overline{x})$, which contradicts Slater's assumption.
\end{dem}

\bigskip

Before we close this section we consider now\ the setting\ of Banach spaces
and we ask whether formula (\ref{generalformula}) is still valid when the
$\limsup$ is taken with respect to the norm topology. As expected, the answer
is affirmative in\ reflexive Banach spaces as the following corollary shows:

\begin{cor}
\label{corref}With the notation of Theorem \emph{\ref{tmain}} we assume that
$(X,\left\Vert ~\right\Vert )$ is a reflexive Banach space. Then\ we have\
\begin{equation}
\mathrm{N}_{[\Phi\leq\Phi(\bar{x})]}^{\delta}(\bar{x})=\left\Vert ~\right\Vert
_{\ast}\text{-}\limsup_{\substack{\sum\nolimits_{i\in\overline{1,k}}\mu
_{i}(\Phi(\bar{x})-\Phi_{t_{i}}(\bar{x})+\varepsilon_{i})\rightarrow
\delta\\\mu_{i}\geq0,\text{ }\lambda_{i}\uparrow\lambda,\text{ }t_{i}\in
T,\text{ }k\in\mathbb{N}}}\sum\nolimits_{i\in\overline{1,k}}\mu_{i}%
\partial_{\varepsilon_{i}}\Phi_{t_{i}}(\bar{x}), \label{nn}%
\end{equation}
where $\left\Vert ~\right\Vert _{\ast}$ denotes the dual norm in $X^{\ast}.$
\end{cor}

\begin{dem}
We consider the pair $(X,X^{\ast})$ with $X^{\ast}$ being the topological dual
space of $X$ associated with the norm topology. Assume, for simplicity, that
$\Phi\equiv\Phi_{t}$ for all $t\in T.$ So, according to Corollary \ref{cor0},
we only need to check that
\begin{equation}
\limsup_{\substack{\mu(\lambda-\Phi(\bar{x})+\varepsilon)\rightarrow
\delta\\\mu\geq0}}\mu\partial_{\varepsilon}\Phi(\overline{x})\subset\left\Vert
~\right\Vert _{\ast}\text{-}\limsup_{\substack{\mu(\lambda-\Phi(\bar
{x})+\varepsilon)\rightarrow\delta\\\mu\geq0}}\mu\partial_{\varepsilon}%
\Phi(\overline{x}). \label{inctopo}%
\end{equation}
Also, by\ Lemma \ref{lem2} given in the Appendix, it suffices to suppose that
$\bar{x}\in\operatorname{argmin}\Phi.$ Take $x^{\ast}$ in the left hand-side
and let nets $(\varepsilon_{\alpha})_{\alpha},$ $(\mu_{\alpha})_{\alpha}$ and
$(x_{\alpha})_{\alpha}$ be such that $\varepsilon_{\alpha},$ $\mu_{\alpha}>0,$
$x_{\alpha}^{\ast}\in\partial_{\varepsilon_{\alpha}}\Phi(\overline{x}),$
$\mu_{\alpha}\varepsilon_{\alpha}\rightarrow\delta,$ and $\mu_{\alpha
}x_{\alpha}^{\ast}\rightharpoonup x^{\ast}.$ Fix $\eta>0$ and let
$\mathcal{A}$ be the set of elements $\alpha$ such that $\mu_{\alpha
}\varepsilon_{\alpha}\leq\delta+\eta.$ Then $x^{\ast}\in\overline
{\operatorname*{co}}\left\{  \mu_{\alpha}x_{\alpha}^{\ast}\right\}
_{\alpha\in\mathcal{A}}=\overline{\operatorname*{co}}^{\Vert~\Vert_{\ast}%
}\left\{  \mu_{\alpha}x_{\alpha}^{\ast}\right\}  _{\alpha\in\mathcal{A}},$ due
to\ the reflexivity of $X,$ and\ there exists $x_{\eta}^{\ast}:=\sum
_{i=1}^{n_{\eta}}\lambda_{i}^{\eta}\mu_{\alpha_{i}}x_{\alpha_{i}}^{\ast}$
$(n_{\eta}\geq1$ and $(\lambda_{i}^{\eta})\in\Delta_{n_{\eta}})$ such that
$\Vert x^{\ast}-x_{\eta}^{\ast}\Vert\leq\eta.$ Since $\bar{x}\in
\operatorname{argmin}\Phi,$ for each $i\in\overline{1,n_{\eta}}$ we have,
setting $\mu_{\eta}^{\prime}:=\max_{j\in\overline{1,n_{\eta}}}\mu_{\alpha_{j}%
}$ $(>0)$ and $\varepsilon_{\eta}^{\prime}:=\frac{\max_{k\in\overline
{1,n_{\eta}}}\varepsilon_{\alpha_{k}}\mu_{\alpha_{k}}}{\max_{j\in
\overline{1,n_{\eta}}}\mu_{\alpha_{j}}},$
\[
\left(  \frac{\mu_{\alpha_{i}}}{\mu_{\eta}^{\prime}}\right)  x_{\alpha_{i}%
}^{\ast}\in\partial_{\frac{\varepsilon_{\alpha_{i}}\mu_{\alpha_{i}}}{\mu
_{\eta}^{\prime}}}(\frac{\mu_{\alpha_{i}}}{\mu_{\eta}^{\prime}}\Phi
)(\overline{x})\subset\partial_{\frac{\varepsilon_{\alpha_{i}}\mu_{\alpha_{i}%
}}{\mu_{\eta}^{\prime}}}\Phi(\overline{x})\subset\partial_{\frac{\max
_{k\in\overline{1,n_{\eta}}}\varepsilon_{\alpha_{k}}\mu_{\alpha_{k}}}%
{\mu_{\eta}^{\prime}}}\Phi(\overline{x})=\partial_{\varepsilon_{\eta}^{\prime
}}\Phi(\overline{x}),
\]
which shows that $x_{\eta}^{\ast}\in\mu_{\eta}^{\prime}\partial_{\varepsilon
_{\eta}^{\prime}}\Phi(\overline{x}).$ But $\mu_{\eta}^{\prime}\varepsilon
_{\eta}^{\prime}=\max_{k\in\overline{1,n_{\eta}}}\mu_{\alpha_{k}}%
\varepsilon_{\alpha_{k}}\leq\delta+\eta$ and, so, by choosing a convergent
subnet of $(\mu_{\eta}^{\prime}\varepsilon_{\eta}^{\prime})_{\eta}$ we arrive
at
\[
x^{\ast}\in\left\Vert ~\right\Vert _{\ast}\text{-}\limsup_{\mu\varepsilon
\rightarrow\delta,\text{ }\mu\geq0}\mu\partial_{\varepsilon}\Phi(\overline
{x}).
\]

\end{dem}

\section{Subdifferential approach\label{sect4}}

As in the previous sections, here we also work with a\ dual pair $(X,Y)$ of
lcs $X$ and $Y$, which are\ endowed with compatible topologies. Our aim is to
characterize the normal cone to sublevel sets by using exclusively the exact
subdifferential of the nominal function.

We consider in this section some restrictions either on the underlying space,
or on the nominal function. We proceed in this way because of\ the existence
in every non-complete normed space of convex proper lsc functions (e.g.,
\cite{brondsted1965subdifferentiability}), which have\ empty subdifferential
mapping everywhere. For simplicity of the presentation, we only study\ the
normal cone to the sublevel set $[\Phi\leq\Phi(\bar{x})].$

In what follows, $\Phi:X\rightarrow\mathbb{R}\cup\left\{  +\infty\right\}  $
is a proper lsc convex function defined on the lcs $X$. First, we recall the
definition of epi-pointed functions.

\begin{defi}
Function $\Phi$ is said to be epi-pointed if its conjugate $\Phi^{\ast}$ is
finite and Mackey-continuous at least at some point.
\end{defi}

\bigskip

The following lemma gathers some useful properties of epi-pointed functions,
which can be found in \cite[Lemma 2.1.6 in and Theorem 4.2]%
{correa2016bronsted}.

\begin{lem}
\label{lema}Assume that function $\Phi$ is epi-pointed. Then the following
assertions hold\emph{:}

\emph{(i)} For every $\varepsilon>0$ we have
\[
\overline{\partial_{\varepsilon}\Phi(x)\cap\operatorname*{int}%
(\operatorname*{dom}\Phi^{\ast})}=\partial_{\varepsilon}\Phi(x).
\]

\emph{(ii)} Assume that $x_{0}^{\ast}\in\partial_{\varepsilon}\Phi(x_{0}%
)\cap\operatorname*{int}(\operatorname*{dom}\Phi^{\ast})$, for some
$\varepsilon\geq0$ and $x_{0}\in X.$ Then for every $\beta\geq0,$ every
continuous seminorm $p$ in $X$, and every $\lambda>0,$ there are
$x_{\varepsilon}\in X,$ $y_{\varepsilon}^{\ast}\in\left[  p\leq1\right]
^{\circ}$ and $\lambda_{\varepsilon}\in\lbrack-1,1]$ such that\emph{:}%
\[
p(x_{0}-x_{\varepsilon})+\beta\left\vert \left\langle x_{0}^{\ast}%
,x_{0}-x_{\varepsilon}\right\rangle \right\vert \leq\lambda,\text{ }\left\vert
\left\langle x_{\varepsilon}^{\ast},x_{0}-x_{\varepsilon}\right\rangle
\right\vert \leq\varepsilon+\frac{\lambda}{\beta},\text{ }|\Phi(x_{0}%
)-\Phi(x_{\varepsilon})|\leq\varepsilon+\frac{\lambda}{\beta},\text{ and }%
\]%
\[
x_{\varepsilon}^{\ast}:=x_{0}^{\ast}+\frac{\varepsilon}{\lambda}%
(y_{\varepsilon}^{\ast}+\beta\lambda_{\varepsilon}x_{0}^{\ast})\in\partial
\Phi(x_{\varepsilon})\cap\partial_{2\varepsilon}\Phi(x_{0})
\]
(with the convention that $\frac{1}{+\infty}=0$).
\end{lem}

\bigskip

\begin{lem}
\label{teoepi}Assume that function $\Phi$\ is epi-pointed. Then for every
$\bar{x}\in\operatorname*{dom}\Phi$
\[
\mathrm{N}_{[\Phi\leq\Phi(\bar{x})]}(\bar{x})=\limsup_{\substack{x\rightarrow
\bar{x},\text{ }\mu\geq0\\\mu(\Phi(x)-\Phi(\bar{x}))\rightarrow0\\\mu
\left\langle \cdot,x-\bar{x}\right\rangle \rightarrow0}}\mu\partial\Phi(x);
\]
that is, each $x^{\ast}\in\mathrm{N}_{[\Phi\leq\Phi(\bar{x})]}(\bar{x})$ is
the (weak*-)limit of a net $(\mu_{i}x_{i}^{\ast}),$\ $\mu_{i}\geq0$ and
$x_{i}^{\ast}\in\partial\Phi(x_{i}),$ with $x_{i}\rightarrow\bar{x},$ $\mu
_{i}(\Phi(x_{i})-\Phi(\bar{x}))\rightarrow0,$ and $\mu\left\langle x_{i}%
^{\ast},x-\bar{x}\right\rangle \rightarrow0.$
\end{lem}

\begin{dem}
The inclusion \textquotedblleft$\supset$\textquotedblright\ being direct, we
are going to prove the inclusion \textquotedblleft$\subset$\textquotedblright.
Given $\xi\in\mathrm{N}_{[\Phi\leq\Phi(\bar{x})]}(\bar{x}),$ we fix\ two open
neighborhoods $U\in\mathcal{N}_{X}(\theta)$ and $V\in\mathcal{N}_{X^{\ast}%
}(\theta),$ and choose $W\in\mathcal{N}_{X^{\ast}}(\theta)$ such that
$W+W\subset V$. According\ to Theorem \ref{tmain}, for each $\varepsilon>0$
there exist $\mu>0$ and $x_{0}^{\ast}\in\partial_{\varepsilon}(\mu\Phi
)(\bar{x})$ such that $x_{0}^{\ast}\in\xi+W$. Let\ $p_{U}:X\rightarrow
\mathbb{R}_{+}$ denote the Minkowski functional defined as\
\[
p_{U}(x):=\inf\left\{  t>0\mid x\in tU\right\}  ,
\]
which is a continuous seminorm on $X.$ Similarly, we define the
weak*-continuous seminorm $p_{V}:X^{\ast}\rightarrow\mathbb{R}_{+}.$ By Lemma
\ref{lema}, we can assume that $x_{0}^{\ast}\in\partial_{\varepsilon}(\mu
\Phi)(\bar{x})\cap\operatorname*{int}(\operatorname*{dom}(\mu\Phi)^{\ast})$,
so that from Lemma\ \ref{lema}, applied with $\lambda=\sqrt{\varepsilon}$ and
$\beta=1$, there exist $x_{\varepsilon}\in X$ and $x_{\varepsilon}^{\ast}%
\in\partial\Phi(x_{\varepsilon})$ such that
\begin{align}
p_{U}(\bar{x}-x_{\varepsilon})  &  \leq\sqrt{\varepsilon},\label{ex}\\
p_{V}(\mu x_{\varepsilon}^{\ast}-(1+\sqrt{\varepsilon})x_{0}^{\ast})  &
\leq\sqrt{\varepsilon}\sup_{y^{\ast}\in U^{\circ}}p_{V}(y^{\ast}%
),\label{exstar}\\
\left\vert \mu\left\langle x_{\varepsilon}^{\ast},x_{\varepsilon}-\bar
{x}\right\rangle \right\vert  &  \leq\varepsilon+\sqrt{\varepsilon}%
,\label{c1}\\
|\mu\Phi(x_{\varepsilon})-\mu\Phi(\bar{x})|  &  \leq\varepsilon+\sqrt
{\varepsilon}. \label{c2}%
\end{align}
If $\varepsilon>0$ is small enough such that $\sqrt{\varepsilon}<1,$ then
condition (\ref{ex}) ensures that $x_{\varepsilon}\in\bar{x}+U.$ Moreover,
since $x_{0}^{\ast}-\xi\in W$ we infer that $2p_{V}(x_{0}^{\ast}-\xi)\leq
p_{W}(x_{0}^{\ast}-\xi)\leq1.$ But $U^{\circ}$ is weak*-compact, by the
Banach-Alaoglu-Bourbaki theorem, and $p_{V}$ is weak*-continuous, and so the
supremum in (\ref{exstar}) is finite and we may assume that $p_{V}(x_{0}%
^{\ast}-\mu x_{\varepsilon}^{\ast})\leq\sqrt{\varepsilon}p_{V}(x_{0}^{\ast
})+\sqrt{\varepsilon}<\frac{1}{2}.$ Consequently,
\[
p_{V}(\mu x_{\varepsilon}^{\ast}-\xi)\leq p_{V}(x_{0}^{\ast}-\mu
x_{\varepsilon}^{\ast})+p_{V}(x_{0}^{\ast}-\xi)<1,
\]
and we get $\mu x_{\varepsilon}^{\ast}-\xi\in V$. Finally, with the use of
conditions (\ref{c1}) and (\ref{c2}) we can easily see that $\xi$ belongs
to\ the right-hand side of the desired formula.
\end{dem}

\bigskip

\begin{cor}
\label{corsin}With the notation of Lemma \emph{\ref{teoepi}} we have
\[
\mathrm{N}_{[\Phi\leq\Phi(\bar{x})]}(\bar{x})=\limsup_{\substack{\mu
(\Phi(x)-\Phi(\bar{x}))\rightarrow0\\\mu\left\langle \cdot,x-\bar
{x}\right\rangle \rightarrow0,\text{ }\mu\geq0}}\mu\partial\Phi(x).
\]

\end{cor}

\begin{dem}
Take $x^{\ast}$ in the right-hand side, so that $x^{\ast}=\lim_{i}\mu_{i}%
x_{i}^{\ast}$ for some $\mu_{i}\geq0$ and $(x_{i},x_{i}^{\ast})\in\partial
\Phi$ such that $\mu_{i}(\Phi(x_{i})-\Phi(\bar{x}))\rightarrow0$ and $\mu
_{i}\left\langle \cdot,x_{i}-\bar{x}\right\rangle \rightarrow0.$ Then, given
$y\in\lbrack\Phi\leq\Phi(\bar{x})]$ and $\varepsilon>0,$ for each $i$
(w.l.o.g.)\ we obtain
\begin{align*}
\left\langle \mu_{i}x_{i}^{\ast},y-\bar{x}\right\rangle  &  =\left\langle
\mu_{i}x_{i}^{\ast},y-x_{i}\right\rangle +\left\langle \mu_{i}x_{i}^{\ast
},x_{i}-\bar{x}\right\rangle \\
&  \leq\mu_{i}(\Phi(y)-\Phi(\bar{x}))+\mu_{i}(\Phi(\bar{x})-\Phi
(x_{i}))+\left\langle \mu_{i}x_{i}^{\ast},x_{i}-\bar{x}\right\rangle
\leq\varepsilon.
\end{align*}
This entails that $x^{\ast}\in\mathrm{N}_{[\Phi\leq\Phi(\bar{x})]}(\bar{x}),$
whereas Lemma \ref{teoepi} gives us
\[
\mathrm{N}_{[\Phi\leq\Phi(\bar{x})]}(\bar{x})=\limsup_{\substack{x\rightarrow
\bar{x},\text{ }\mu\geq0\\\mu(\Phi(x)-\Phi(\bar{x}))\rightarrow0\\\mu
\left\langle \cdot,x-\bar{x}\right\rangle \rightarrow0}}\mu\partial
\Phi(x)\subset\limsup_{\substack{\mu(\Phi(x)-\Phi(\bar{x}))\rightarrow
0\\\mu\left\langle \cdot,x-\bar{x}\right\rangle \rightarrow0,\text{ }\mu\geq
0}}\mu\partial\Phi(x)\subset\mathrm{N}_{[\Phi\leq\Phi(\bar{x})]}(\bar{x}).
\]

\end{dem}

Up to a little modification in the way that one perturbs the nominal point
$\bar{x},$ Lemma \ref{teoepi} above is still valid if the epi-pointedness
condition is weakened to the continuity of the conjugate relative to its
domain. We obtain then the following result, which gives both primal and dual
symmetric condition ensuring the same characterization of the normal cone. One
can understand that these conditions preclude the subdifferential mapping from
being empty everywhere.

The notation $x\rightarrow_{A}\bar{x},$ $A\subset X^{\ast},$ used below,
refers to the convergence of the corresponding equivalence classes in the
quotient space $X/A^{\bot}.$

\begin{thm}
\label{rii}Assume that either $\Phi_{\mid\operatorname*{aff}%
(\operatorname*{dom}\Phi)}$ is finite and continuous in $\operatorname*{ri}%
(\operatorname*{dom}\Phi)$ $(\neq\emptyset)$ or $\Phi_{\mid\operatorname*{aff}%
(\operatorname*{dom}\Phi^{\ast})}^{\ast}$ is finite and (Mackey-)continuous in
$\operatorname*{ri}(\operatorname*{dom}\Phi^{\ast})\ (\neq\emptyset).$ Then
for every $\bar{x}\in\operatorname*{dom}\Phi$
\[
\mathrm{N}_{[\Phi\leq\Phi(\bar{x})]}(\bar{x})=\limsup_{\substack{x\rightarrow
_{\operatorname*{dom}\Phi^{\ast}}\bar{x},\text{ }\mu\geq0\\\mu(\Phi
(x)-\Phi(\bar{x}))\rightarrow0\\\mu\left\langle \cdot,x-\bar{x}\right\rangle
\rightarrow0}}\mu\partial\Phi(x)=\limsup_{\substack{\mu(\Phi(x)-\Phi(\bar
{x}))\rightarrow0\\\mu\left\langle \cdot,x-\bar{x}\right\rangle \rightarrow
0,\text{ }\mu\geq0}}\mu\partial\Phi(x).
\]

\end{thm}

\begin{dem}
We assume first that $\Phi_{\mid\operatorname*{aff}(\operatorname*{dom}%
\Phi^{\ast})}^{\ast}$ is finite and (Mackey-)continuous in the set
$\operatorname*{ri}(\operatorname*{dom}\Phi^{\ast})\ (\neq\emptyset).$ In view
of Lemma \ref{lem2} we may assume that $\theta\in\partial\Phi(\bar{x});$
otherwise, we are obviously done. Hence, $\theta\in\operatorname*{dom}%
\Phi^{\ast}\subset\operatorname*{aff}(\operatorname*{dom}\Phi^{\ast})$ so that
$Y^{\ast}:=\operatorname*{aff}(\operatorname*{dom}\Phi^{\ast})$ is a closed
subspace of $X^{\ast}.$ Denote $Y:=X/\operatorname*{aff}(\operatorname*{dom}%
\Phi^{\ast})^{\bot}$ (the quotient space) so that the pair $(Y,Y^{\ast})$
becomes a dual pair when endowed with the quotient and trace topologies of $X$
and $X^{\ast},$ respectively. Consider the function $\widetilde{\Phi}$ defined
on $Y$ as $\widetilde{\Phi}(\tilde{x})=\Phi(x)$ where $x$ is in the equivalent
class of $\tilde{x}\in Y.$ It can be easily checked that $\widetilde{\Phi}%
\in\Gamma_{0}(Y)$ and that $(\widetilde{\Phi})^{\ast}(y^{\ast})=\Phi^{\ast
}(y^{\ast})$ for all $y^{\ast}\in\operatorname*{dom}\Phi^{\ast}%
(=\operatorname*{dom}(\widetilde{\Phi})^{\ast});$ hence, $\widetilde{\Phi}$ is
epi-pointed too. Take $x^{\ast}\in\mathrm{N}_{[\Phi\leq\Phi(\bar{x})]}(\bar
{x})$ so that, by Lemma \ref{teoepi}, $x^{\ast}\in\operatorname*{cl}%
(\mathbb{R}_{+}(\operatorname*{range}\partial\Phi))\subset\operatorname*{cl}%
(\mathbb{R}_{+}(\operatorname*{dom}\Phi^{\ast}))\subset Y^{\ast}$ and, hence,
$x^{\ast}\in\mathrm{N}_{[\widetilde{\Phi}\leq\widetilde{\Phi}(\widetilde{\bar
{x}})]}(\widetilde{\bar{x}}).$ By applying again Lemma \ref{teoepi}, and also
Corollary \ref{corsin}, we get
\[
x^{\ast}\in\limsup_{\substack{\tilde{x}\rightarrow\widetilde{\bar{x}},\text{
}\mu\geq0\\\mu(\widetilde{\Phi}(\tilde{x})-\widetilde{\Phi}(\widetilde{\bar
{x}}))\rightarrow0\\\mu\left\langle \cdot,\tilde{x}-\widetilde{\bar{x}%
}\right\rangle \rightarrow0}}\mu\partial\widetilde{\Phi}(\tilde{x}%
)=\limsup_{\substack{x\rightarrow_{\operatorname*{dom}\Phi^{\ast}}\bar
{x},\text{ }\mu\geq0\\\mu(\Phi(x)-\Phi(\bar{x}))\rightarrow0\\\mu\left\langle
\cdot,x-\bar{x}\right\rangle \rightarrow0}}\mu\partial\Phi(x)\subset
\limsup_{\substack{\mu(\Phi(x)-\Phi(\bar{x}))\rightarrow0\\\mu\left\langle
\cdot,x-\bar{x}\right\rangle \rightarrow0,\text{ }\mu\geq0}}\mu\partial
\Phi(x),
\]
and we conclude.

Assume now that $\Phi_{\mid\operatorname*{aff}(\operatorname*{dom}\Phi)}$ is
finite and continuous in $\operatorname*{ri}(\operatorname*{dom}\Phi).$ Let
$L$ be a finite-dimensional subspace of $X$ which contains $\bar{x}$ and
denote by\ $\Phi_{L}$ the restriction of $\Phi$ to $L;$ hence, $\Phi_{L}%
^{\ast}\in\Gamma_{0}(L^{\ast})$ $(L^{\ast}$ being the dual space of $L)$ and
$\operatorname*{ri}(\operatorname*{dom}\Phi_{L}^{\ast})\neq\emptyset.$ From
the first part of the proof applied in the dual pair $(L,L^{\ast})$ we get
\[
\mathrm{N}_{[\Phi_{L}\leq\Phi(\bar{x})]}(\bar{x})=\limsup
_{\substack{x\rightarrow_{\operatorname*{dom}\Phi_{L}^{\ast}}\bar{x},\text{
}\mu\geq0\\\mu(\Phi_{L}(x)-\Phi(\bar{x}))\rightarrow0\\\mu\left\langle
\cdot,x-\bar{x}\right\rangle \rightarrow0}}\mu\partial\Phi_{L}(x).
\]
We pick an\ $x^{\ast}\in\mathrm{N}_{[\Phi\leq\Phi(\bar{x})]}(\bar{x});$ hence,
$x_{\mid L}^{\ast}\in\mathrm{N}_{[\Phi_{L}\leq\Phi(\bar{x})]}(\bar{x}).$ From
the last relation above we find $x_{i}\in L$ $(\subset X),$ $\mu_{i}\geq0$ and
$\hat{x}_{i}^{\ast}\in\partial\Phi_{L}(x_{i})$ such that $x_{\mid L}^{\ast
}=w^{\ast}$-$\lim_{i}\mu_{i}\hat{x}_{i}^{\ast},$ $x_{i}\rightarrow
_{\operatorname*{dom}\Phi_{L}^{\ast}}\bar{x},$ and
\[
\mu_{i}\left\langle \hat{x}_{i}^{\ast},x_{i}-\bar{x}\right\rangle
\rightarrow0,\ \mu_{i}(\Phi(x_{i})-\Phi(\bar{x}))=\mu_{i}(\Phi_{L}(x_{i}%
)-\Phi(\bar{x}))\rightarrow0.
\]
By Hahn-Banach we extend $\hat{x}_{i}^{\ast}$ to $x_{i}^{\ast}\in X^{\ast}$
such that $x_{i\mid L}^{\ast}\equiv\hat{x}_{i}^{\ast}$ and the net $(\mu
_{i}x_{i}^{\ast})_{i}$ remains bounded, so, weak* convergent (w.lo.g.). It
follows that\ $x_{i}^{\ast}\in\partial(\Phi+\mathrm{I}_{L})(x_{i}%
)=\partial\Phi(x_{i})+L^{\perp}$ (e.g., \cite{HanLopCor16}), $\mu
_{i}\left\langle x_{i}^{\ast},x_{i}-\bar{x}\right\rangle \rightarrow0$ and
\[
x^{\ast}\in(w^{\ast}\text{-}\lim_{i}\mu_{i}x_{i}^{\ast})+L^{\perp}.
\]
Finally, since $(\operatorname*{dom}\Phi^{\ast})_{\mid L}:=\{u_{\mid L}^{\ast
}\mid u^{\ast}\in\operatorname*{dom}\Phi^{\ast}\}\subset\operatorname*{dom}%
\Phi_{L}^{\ast},$ the finite-dimensionality of $L$ implies that $x_{i}%
\rightarrow_{\operatorname*{dom}\Phi^{\ast}}\bar{x}.$ In other words, due to
the arbitrariness of $L,$ we conclude that $x^{\ast}\in\limsup
_{\substack{x\rightarrow_{\operatorname*{dom}\Phi^{\ast}}\bar{x},\text{ }%
\mu\geq0\\\mu(\Phi(x)-\Phi(\bar{x}))\rightarrow0\\\mu\left\langle \cdot
,x-\bar{x}\right\rangle \rightarrow0}}\mu\partial\Phi(x),$ which yields the
left inclusion \textquotedblleft$\subset$\textquotedblright\ of the required
statement. The other inclusions are easy and have been proved in previous
opportunities (see Corollary \ref{corsin}).
\end{dem}

\bigskip

It is worth observing that in the Banach setting any convex function can be
made epi-pointed via a penalization with the indicator function of a bounded
set. Then we obtain the following result, given originally in
\cite{cabot2014sequential}.

\begin{cor}
\label{ct}Assume that $X$ is Banach. Then for every $\bar{x}\in
\operatorname*{dom}\Phi$ we have that
\[
\mathrm{N}_{[\Phi\leq\Phi(\bar{x})]}(\bar{x})=\limsup_{\substack{x\rightarrow
\bar{x},\text{ }\mu\geq0\\\mu(\Phi(x)-\Phi(\bar{x}))\rightarrow0\\\mu
\left\langle \cdot,x-\bar{x}\right\rangle \rightarrow0}}\mu\partial\Phi(x),
\]
where the limit is taken with respect to the weak*-topology.
\end{cor}

\begin{dem}
We consider the dual pair $(X,X^{\ast})$, where $X^{\ast}$ is the topological
dual of $X^{\ast}$ endowed with the weak*-topology. Fix a non-zero
element\ $\xi\in\mathrm{N}_{[\Phi\leq\Phi(\bar{x})]}(\bar{x})$ and pick a
$\theta$-neighborhood $V$ (with respect to the Mackey-topology in $X^{\ast})$
together with an $\eta>0.$ Next, we choose a weakly compact (and symmetric)
convex set $K\subset X$ such that $8K^{\circ}\subset V$. We set
$\widetilde{\Phi}:=\Phi+\mathrm{I}_{\bar{x}+K}$ $(\in\Gamma_{0}(X)).$ It is
easy to verify\ that the conjugate function $\widetilde{\Phi}^{\ast}$ is
finite and bounded from above on $x_{0}^{\ast}+K^{\circ}$ for some
$x_{0}^{\ast}\in\operatorname*{dom}\Phi^{\ast},$ and, so, is
(Mackey-)continuous on $x_{0}^{\ast}+K.$ It follows that $\widetilde{\Phi}$ is
epi-pointed, by definition, so that Lemma \ref{teoepi} applies and yields\
\[
\xi\in\mathrm{N}_{[\Phi\leq\Phi(\bar{x})]}(\bar{x})\subset\mathrm{N}%
_{[\Phi\leq\Phi(\bar{x})]\cap(\bar{x}+K)}(\bar{x})=\mathrm{N}%
_{[\widetilde{\Phi}\leq\Phi(\bar{x})]}(\bar{x})=\limsup
_{\substack{x\rightarrow\bar{x},\text{ }\mu\geq0\\\mu(\Phi(x)-\Phi(\bar
{x}))\rightarrow0\\\mu\left\langle \cdot,x-\bar{x}\right\rangle \rightarrow
0}}\partial(\mu\widetilde{\Phi})(x).
\]
Thus,\ for any\ $\theta$-neighborhood $U\in\mathcal{N}_{X}(\theta),$ there
exist $\mu_{i}>0,$ $\hat{x}_{i}\in\bar{x}+U,$ and $\hat{x}_{i}^{\ast}%
\in\partial(\mu_{i}\widetilde{\Phi})(\hat{x}_{i})$ such that
\begin{equation}
\xi\in\hat{x}_{i}^{\ast}+\eta K^{\circ},\text{ }\left\vert \mu_{i}(\Phi
(\hat{x}_{i})-\Phi(\bar{x}))\right\vert \leq\eta, \label{st8}%
\end{equation}%
\begin{equation}
\left\vert \left\langle \hat{x}_{i}^{\ast},\hat{x}_{i}-\bar{x}\right\rangle
\right\vert \leq\eta; \label{st3}%
\end{equation}
in particular, we have that
\begin{equation}
\hat{x}_{i}\in\bar{x}+U\cap K. \label{st4}%
\end{equation}
On another hand, by the sum rule in \cite{thibault1997sequential},\ for each
$i$ we\ find
\begin{equation}
x_{i}\in\hat{x}_{i}+U,\text{ }x_{i}^{\ast}\in\partial(\mu_{i}\Phi
)(x_{i}),\text{ } \label{st6}%
\end{equation}%
\begin{equation}
y_{i}^{\ast}\in\partial_{\eta}(\mu_{i}\mathrm{I}_{\bar{x}+K})(\hat{x}_{i}%
)=\mu_{i}\mathrm{N}_{\bar{x}+K}^{\eta\mu_{i}^{-1}}(\hat{x}_{i}) \label{st1}%
\end{equation}
such that
\begin{equation}
\hat{x}_{i}^{\ast}\in x_{i}^{\ast}+y_{i}^{\ast}+\eta K^{\circ}, \label{st2}%
\end{equation}%
\begin{equation}
\left\vert \left\langle x_{i}^{\ast},x_{i}-\hat{x}_{i}\right\rangle
\right\vert \leq\eta, \label{st5}%
\end{equation}%
\[
\max\{\left\vert \mu_{i}(\Phi(\hat{x}_{i})-\Phi(x_{i}))\right\vert
,~\left\vert \left\langle x_{i}^{\ast},\hat{x}_{i}-x_{i}\right\rangle
\right\vert \}\leq\eta;
\]
hence, in particular,
\begin{equation}
\{\left\vert \mu_{i}(\Phi(x_{i})-\Phi(\bar{x}))\right\vert ,\left\vert
\left\langle \hat{x}_{i}^{\ast},x_{i}-\bar{x}\right\rangle \right\vert
\}\leq2\eta. \label{st7}%
\end{equation}
Now, for every $i$ and $z\in K$ it holds%
\[%
\begin{array}
[c]{llll}%
\left\langle y_{i}^{\ast},z\right\rangle  & \leq & \eta+\left\langle
y_{i}^{\ast},\hat{x}_{i}-\bar{x}\right\rangle  & \text{by (\ref{st1})}\\
& = & \eta+\left\langle x_{i}^{\ast}+y_{i}^{\ast}-\hat{x}_{i}^{\ast},\hat
{x}_{i}-\bar{x}\right\rangle +\left\langle \hat{x}_{i}^{\ast},\hat{x}_{i}%
-\bar{x}\right\rangle -\left\langle x_{i}^{\ast},\hat{x}_{i}-\bar
{x}\right\rangle  & \\
& \leq & \eta+\eta\mathrm{\sigma}_{K^{\circ}}(\hat{x}_{i}-\bar{x}%
)+\left\langle \hat{x}_{i}^{\ast},\hat{x}_{i}-\bar{x}\right\rangle
-\left\langle x_{i}^{\ast},\hat{x}_{i}-\bar{x}\right\rangle  & \text{by
(\ref{st2})}\\
& \leq & \eta+\eta+\eta-\left\langle x_{i}^{\ast},\hat{x}_{i}-\bar
{x}\right\rangle  & \text{by (\ref{st4}) and (\ref{st3})}\\
& \leq & 3\eta-\left\langle x_{i}^{\ast},\hat{x}_{i}-\bar{x}\right\rangle  &
\text{by (\ref{st1})}\\
& = & 3\eta-\left\langle x_{i}^{\ast},\hat{x}_{i}-x_{i}\right\rangle
+\left\langle x_{i}^{\ast},\bar{x}-x_{i}\right\rangle  & \\
& \leq & 3\eta+\eta+\left\langle x_{i}^{\ast},\bar{x}-x_{i}\right\rangle  &
\text{by (\ref{st5})}\\
& \leq & 4\eta+(\mu_{i}\Phi)(\bar{x})-(\mu_{i}\Phi)(x_{i}) & \text{by
(\ref{st6})}\\
& \leq & 4\eta+2\eta=6\eta & \text{by (\ref{st7}),}%
\end{array}
\]
so that
\begin{equation}
y_{i}^{\ast}\in6\eta K^{\circ}; \label{hl}%
\end{equation}
hence, it also follows from the seventh inequality in the table above that
\begin{equation}
\left\langle x_{i}^{\ast},x_{i}-\bar{x}\right\rangle \leq4\eta+\mathrm{\sigma
}_{K}(y_{i}^{\ast})\leq10\eta. \label{hl2}%
\end{equation}
Consequently, using successively (\ref{st8}), (\ref{st2}), together with the
choice of $V$ in the beginning of the proof, (\ref{hl}) gives us
\begin{equation}
\xi\in\hat{x}_{i}^{\ast}+\eta K^{\circ}\subset x_{i}^{\ast}+y_{i}^{\ast}+2\eta
K^{\circ}\subset x_{i}^{\ast}+8\eta K^{\circ}\subset x_{i}^{\ast}+V.
\label{gu}%
\end{equation}
Finally, since\ $x_{i}^{\ast}\in\partial(\mu_{i}\Phi)(x_{i})$ (recall
(\ref{st6})) and we have that
\[
x_{i}\in\hat{x}_{i}+U\in\bar{x}+U+U,\text{ by (\ref{st6}) and (\ref{st4}),}%
\]%
\[
\left\vert \mu(\Phi(x_{i})-\Phi(\bar{x}))\right\vert \leq2\eta,\text{ by
(\ref{st7}),}%
\]%
\[
-2\eta\leq\mu(\Phi(x_{i})-\Phi(\bar{x}))\leq\left\langle x_{i}^{\ast}%
,x_{i}-\bar{x}\right\rangle \leq10\eta,\text{ by (\ref{hl2}), (\ref{st6}), and
(\ref{st7}),}%
\]
together with (\ref{st7}), relation (\ref{gu}) leads us to $\xi\in
\limsup_{\substack{x\rightarrow\bar{x},\text{ }\mu>0\\\mu(\Phi(x)-\Phi(\bar
{x}))\rightarrow0\\\mu\left\langle \cdot,x-\bar{x}\right\rangle \rightarrow
0}}\mu\partial\Phi(x)+V.$ The desired inclusion follows then by the
arbitrariness of $V.$
\end{dem}

\bigskip

\begin{rem}
\emph{It is possible to obtain Corollary \ref{ct} directly from Lemma
\ref{teoepi} by applying Borwein's version of Brøndsted-Rockafellar's theorem.
Nevertheless, our approach permits to highlight the generality of Lemma
\ref{teoepi} (and Theorem \ref{rii}) in the sense that the epi-pointedness
condition is not so restrictive as it may appear from a first glance. Let us,
for completeness, give the direct proof:}

\emph{Take }$x^{\ast}$\emph{ in }$\mathrm{N}_{[\Phi\leq\Phi(\bar{x})]}(\bar
{x})$\emph{. By Corollary \ref{cor1} for each }$\delta>0$\emph{ we have}
\begin{equation}
x^{\ast}\in\mathrm{N}_{[\Phi\leq\Phi(\bar{x})]}^{\delta}(\bar{x}%
)=\overline{\bigcup_{\mu>0}\partial_{\delta}(\mu\Phi)(\bar{x})}^{w^{\ast}}.
\label{bronineq1}%
\end{equation}
\emph{Thus, for every weak* neighborhood }$W$\emph{ of }$x^{\ast},$\emph{
there are some }$\mu>0$\emph{ and }$x_{0}^{\ast}\in\partial_{\delta}(\mu
\Phi)(\bar{x})\cap W$\emph{. By Br$\not o  $ndsted-Rockafellar's-like Theorem
(\cite{penot1996subdifferential})\ applied to the function }$f:=\mu\Phi
$\emph{, we find }$x_{\delta}\in X$\emph{ and }$y_{\delta}^{\ast}=\mu
x_{\delta}^{\ast}\in\mu\partial\Phi(x_{\delta})=\partial(\mu\Phi)(x_{\delta}%
)$\emph{ such that }%
\begin{align*}
\Vert x_{\delta}-\bar{x}\Vert &  \leq\sqrt{\delta},\\
\Vert x_{\delta}^{\ast}-x_{0}^{\ast}\Vert &  \leq\sqrt{\delta}(1+\sqrt{\delta
}),\\
|\mu\left\langle x_{\delta}^{\ast},x_{\delta}-\bar{x}\right\rangle |  &
\leq\delta+\sqrt{\delta},\\
|\mu\Phi(x_{\delta})-\mu\Phi(\bar{x})|  &  \leq\delta+\sqrt{\delta}.
\end{align*}
\emph{Since }$W$\emph{ is also open in the norm topology }$\tau_{\Vert
\cdot\Vert_{\ast}}$\emph{, by taking }$\delta$\emph{ small enough if
necessary, the second inequality above guarantees that }$x_{\delta}^{\ast}\in
W$\emph{. So, we get }%
\[
x^{\ast}\in w^{\ast}\text{-}\limsup_{\substack{x\rightarrow\bar{x}\\\mu
(\Phi(x)-\Phi(\bar{x}))\rightarrow0\\\mu\left\langle \cdot,x-\bar
{x}\right\rangle \rightarrow0\\\mu\geq0}}\mu\partial\Phi(x).
\]

\end{rem}

\section{Spectral functions\label{sect5}}

In this last section, we apply Theorems \ref{tmain} and \ref{rii} to derive
characterizations of\ the\ normal cone to sublevel sets of a (proper, convex
and lsc) spectral function, by means of its restriction to the range of the
eigenvalues vectors.

\bigskip

Here, we identify $X$ to the Euclidean space of $n\times n$ symmetric matrices
with coefficients in $\mathbb{R},$ $S^{n}(\mathbb{R})$, which is endowed with
the trace inner product $\left\langle X,Y\right\rangle :=\operatorname*{tr}%
(XY).$ We denote by $S_{-}^{n}(\mathbb{R})$ the cone of semi-definite negative matrices.

An extended real-valued matrix function $F:S^{n}(\mathbb{R})\rightarrow
\mathbb{R}\cup\{+\infty\}$ is said to be spectral if it depends only on the
eigenvalues of the corresponding matrix; that is,
\[
F(A)=F(U^{T}AU),
\]
for any $A\in S^{n}(\mathbb{R})$ and any orthogonal matrix $U$ $($i.e.,
$U\in\mathcal{O}^{n}(\mathbb{R}));$ the superscript $T$ denotes the transpose
matrix. We shall denote by $\mathcal{D}$ the set of (continuous) linear
transformations $A_{U}:\mathbb{R}^{n}\rightarrow S^{n}(\mathbb{R}),$
$U\in\mathcal{O}^{n}(\mathbb{R}),$ given by
\[
A_{U}x:=U^{T}\operatorname*{diag}(x)U,
\]
where $\operatorname*{diag}(x)\ $is the diagonal matrix whose main diagonal is
formed by the elements of vector $x.$

If $F\in\Gamma_{0}(S^{n}(\mathbb{R})),$ then (see \cite{lewis}) there is a
symmetric function $f\in\Gamma_{0}(\mathbb{R}^{n})$ such that $F=f\circ
\lambda,$ where $\lambda(A):=(\lambda_{1}(A),\cdots,\lambda_{n}(A))$ is the
spectra of $A$ arranged in a non-increasing order. We know that
\begin{equation}
(f\circ\lambda)^{\ast}=f^{\ast}\circ\lambda, \label{spe}%
\end{equation}
and, consequently, for every $X\in S^{n}(\mathbb{R})$ and $\varepsilon\geq0$
it holds that
\[
\partial_{\varepsilon}(f\circ\lambda)(X)=A_{U}\partial_{\varepsilon}%
f(\lambda(X)),
\]
where $U\in\mathcal{O}^{n}(\mathbb{R})\ $satisfies\ $X=A_{U}\lambda(X).$
Relation (\ref{spe}) also yields the following comparison between the
recession functions of $f$ and $f\circ\lambda.$

\begin{lem}
\label{lems}Given a symmetric function $f\in\Gamma_{0}(\mathbb{R}^{n}),$ we
have $(f\circ\lambda)^{\infty}=f^{\infty}\circ\lambda.$
\end{lem}

\begin{dem}
Since both functions $f$ and $f\circ\lambda$ are proper, convex and lsc,
(\ref{spe}) gives us $(f\circ\lambda)^{\infty}=\mathrm{\sigma}%
_{\operatorname*{dom}(f\circ\lambda)^{\ast}}=\mathrm{\sigma}%
_{\operatorname*{dom}(f^{\ast}\circ\lambda)}$ and $f^{\infty}\circ
\lambda=\mathrm{\sigma}_{\operatorname*{dom}f^{\ast}}\circ\lambda.$ But
\begin{align*}
\operatorname*{dom}(f^{\ast}\circ\lambda)  &  =\{X\in S^{n}(\mathbb{R}%
)\mid\lambda(X)\in\operatorname*{dom}f^{\ast}\}\\
&  =\{X\in A_{U}\operatorname*{dom}f^{\ast},\text{ }U\in\mathcal{O}%
^{n}(\mathbb{R})\text{ st. }A_{U}\lambda(X)=X\}\\
&  =\{A_{U}\operatorname*{dom}f^{\ast}\mid U\in\mathcal{D}\},
\end{align*}
and, so, using Von Neumann's trace inequality (see, e.g., \cite[Theorem
2.1]{lewis}), for every $X\in S^{n}(\mathbb{R})$%
\begin{align*}
(f\circ\lambda)^{\infty}(X)  &  =\mathrm{\sigma}_{A_{U}\operatorname*{dom}%
f^{\ast},\text{ }U\in D}(X)\\
&  =\sup_{y\in\operatorname*{dom}f^{\ast}}\sup_{U\in\mathcal{O}^{n}%
(\mathbb{R})}\left\langle A_{U}y,X\right\rangle \\
&  =\sup_{y\in\operatorname*{dom}f^{\ast}}\left\langle y,\lambda
(X)\right\rangle )=(f^{\infty}\circ\lambda)(X).
\end{align*}

\end{dem}

\bigskip

We have the following result:

\begin{cor}
\label{diag}Given a symmetric function $f\in\Gamma_{0}(\mathbb{R}^{n})$ and
$\bar{X}\in\operatorname*{dom}(f\circ\lambda),$ we choose a matrix
$U\in\mathcal{O}^{n}(\mathbb{R})$ such that $\bar{X}=A_{U}\lambda(\bar{X}).$
Then for every $\delta\geq0$ and $\alpha\in\left]  -\infty,+\infty\right]  $
we have
\begin{align*}
\mathrm{N}_{([f\circ\lambda\leq\alpha]\cap\operatorname*{dom}(f\circ
\lambda))\cup(\bar{X}+[f^{\infty}\circ\lambda\leq0])}^{\delta}(\bar{X})  &
=A_{U}\mathrm{N}_{([f\leq\alpha]\cap\operatorname*{dom}f)\cup(\lambda(\bar
{X})+[f^{\infty}\leq0])}^{\delta}(\lambda(\bar{X}))\\
&  =A_{U}\limsup_{\substack{\mu(\beta-f(\lambda(\bar{X}))+\varepsilon
)\rightarrow\delta\\\mu\geq0,\text{ }\beta\uparrow\alpha}}\mu\partial
_{\varepsilon}f(\lambda(\bar{X})).
\end{align*}
In addition, if $\delta=0$ and $f(\lambda(\bar{X}))=\alpha,$ then we also
have
\[
\mathrm{N}_{[f\circ\lambda\leq f(\lambda(\bar{X}))]}(\bar{X})=A_{U}%
\limsup_{\substack{x\rightarrow\lambda(\bar{X}),\text{ }\mu\geq0\\\mu
(f(x)-f(\lambda(\bar{X})))\rightarrow0\\\mu\left\langle \cdot,x-\lambda
(\bar{X})\right\rangle \rightarrow0}}\mu\partial f(x).
\]

\end{cor}

\begin{dem}
By applying Theorem \ref{tmain} twice we obtain
\begin{align*}
\mathrm{N}_{([f\circ\lambda\leq\alpha]\cap\operatorname*{dom}(f\circ
\lambda))\cup(\bar{X}+[f^{\infty}\circ\lambda\leq0])}^{\delta}(\bar{X})  &
=\limsup_{\substack{\mu(\beta-f(\lambda(\bar{X}))+\varepsilon)\rightarrow
\delta\\\mu\geq0,\beta\uparrow\alpha}}\mu\partial_{\varepsilon}(f\circ
\lambda)(\bar{X})\\
&  =\limsup_{\substack{\mu(\beta-f(\lambda(\bar{X}))+\varepsilon
)\rightarrow\delta\\\mu\geq0,\beta\uparrow\alpha}}\mu A_{U}\partial
_{\varepsilon}f(\lambda(\bar{X}))\\
&  =A_{U}\limsup_{\substack{\mu(\beta-f(\lambda(\bar{X}))+\varepsilon
)\rightarrow\delta\\\mu\geq0,\beta\uparrow\alpha}}\mu\partial_{\varepsilon
}f(\lambda(\bar{X}))\\
&  =A_{U}\mathrm{N}_{([f\leq\alpha]\cap\operatorname*{dom}f)\cup(\lambda
(\bar{X})+[f^{\infty}\leq0])}^{\delta}(\lambda(\bar{X})),
\end{align*}
which yields the first statement of the corollary. The last statement follows
from Corollary \ref{ct} in a similar way.
\end{dem}

\bigskip

In particular, when $\alpha=f(\lambda(\bar{X}))\in\mathbb{R}$ and $\delta=0,$
the previous characterization gives us
\[
\mathrm{N}_{[f\circ\lambda\leq f(\lambda(\bar{X}))]}(\bar{X})=A_{U}%
\mathrm{N}_{[f\leq f(\lambda(\bar{X}))]}(\lambda(\bar{X}))=A_{U}\limsup
_{\mu\varepsilon\rightarrow0,\text{ }\mu\geq0}\mu\partial_{\varepsilon
}f(\lambda(\bar{X})),
\]
where $U\in\mathcal{O}^{n}(\mathbb{R})$ is such that $\bar{X}=A_{U}%
\lambda(\bar{X}).$ The first equality can be easily obtained from a more
general result given in \cite[Corollary 5.3]{lewisin}, by observing that the
set $C:=[f\circ\lambda\leq f(\lambda(\bar{X}))]$ is an invariant set (see
\cite[Corollary 5.3]{lewisin}). We can proceed similarly when $[f\circ
\lambda\leq\alpha]=\emptyset,$ since that $\mathrm{N}_{([f\circ\lambda
\leq\alpha]\cap\operatorname*{dom}(f\circ\lambda))\cup(\bar{X}+[f^{\infty
}\circ\lambda\leq0])}(\bar{X})=\mathrm{N}_{[f^{\infty}\circ\lambda\leq0]}(0)$
and the set $[f^{\infty}\circ\lambda\leq0]$ is invariant too. However, this
argument can not be used, at least in an immediate way, to get the first
equality in Corollary \ref{diag} when the set $[f\circ\lambda\leq\alpha]$ is
not empty and $f(\lambda(\bar{X}))>\alpha,$ in which case the set
$([f\circ\lambda\leq\alpha]\cap\operatorname*{dom}(f\circ\lambda))\cup(\bar
{X}+[f^{\infty}\circ\lambda\leq0])$ is not necessarily invariant.

\bigskip

We consider in the following corollary the special and typical example of the
largest eigenvalue function, $\lambda_{\max}(X):=\max_{i\in\overline{1,n}%
}\lambda_{i}(X).$

\begin{cor}
Given an $\bar{X}\in S^{n}(\mathbb{R})$ we choose a $U\in\mathcal{O}%
^{n}(\mathbb{R})$ such that $\bar{X}=A_{U}\lambda(\bar{X}).$ Then for every
$\delta\geq0$ and $\alpha\in\mathbb{R}$
\[
\mathrm{N}_{[\lambda_{\max}\leq\alpha]\cup(\bar{X}+S_{-}^{n}(\mathbb{R}%
))}^{\delta}(\bar{X})=A_{U}\limsup_{\substack{\sum\nolimits_{i\in
\overline{1,n}}\mu_{i}(\alpha-\lambda_{\max}(\bar{X})+\varepsilon
)\rightarrow\delta,\text{ }\mu_{i}\geq0\\\sum\limits_{i=1}^{n}\mu_{i}%
(\lambda_{i}(\bar{X})-\lambda_{\max}(\bar{X})-\varepsilon)\geq0}}(\mu
_{1}\lambda_{1}(\bar{X}),\cdots,\mu_{n}\lambda_{n}(\bar{X}))^{T}.
\]

\end{cor}

\begin{dem}
Since $\lambda_{\max}=f\circ\lambda,$ with $f(x_{1},\cdots,x_{n}):=\max
\{x_{1},\cdots,x_{n}\}$ $(\in\Gamma_{0}(\mathbb{R}^{n})),$ Lemma \ref{lems}
ensures that $\lambda_{\max}(X)^{\infty}=f^{\infty}\circ\lambda=f\circ
\lambda=\lambda_{\max}(X).$ Then $[f^{\infty}\circ\lambda\leq0]=S_{-}%
^{n}(\mathbb{R})$ and so, according to Corollary \ref{diag},
\begin{align*}
\mathrm{N}_{[\lambda_{\max}\leq\alpha]\cup(\bar{X}+[f^{\infty}\circ\lambda
\leq0])}^{\delta}(\bar{X})  &  =\mathrm{N}_{[\lambda_{\max}\leq\alpha
]\cup(\bar{X}+S_{-}^{n}(\mathbb{R}))}^{\delta}(\bar{X})\\
&  =A_{U}\limsup_{\mu(\alpha-f(\lambda(\bar{X}))+\varepsilon)\rightarrow
\delta,\text{ }\mu\geq0}\mu\partial_{\varepsilon}f(\lambda(\bar{X})).
\end{align*}
Thus, Corollary \ref{finitecase} leads us to
\begin{align*}
\mathrm{N}_{[\lambda_{\max}\leq\alpha]\cup(\bar{X}+[f^{\infty}\circ\lambda
\leq0])}^{\delta}(\bar{X})  &  =A_{U}\limsup_{\substack{\mu(\alpha
-\lambda_{\max}(\bar{X})+\varepsilon)\rightarrow\delta\\(\gamma_{i})\in
\Delta_{n},\text{ },\text{ }\mu\geq0\\\sum\limits_{i=1}^{n}\gamma_{i}%
\lambda_{i}(\bar{X})\geq\lambda_{\max}(\bar{X})-\varepsilon}}(\mu\gamma
_{1}\lambda_{1}(\bar{X}),\cdots,\mu\gamma_{n}\lambda_{n}(\bar{X}))^{T}\\
&  =A_{U}\limsup_{\substack{\sum\nolimits_{i\in\overline{1,n}}\mu_{i}%
(\alpha-\lambda_{\max}(\bar{X})+\varepsilon)\rightarrow\delta,\mu_{i}%
\geq0\\\sum\limits_{i=1}^{n}\mu_{i}(\lambda_{i}(\bar{X})-\lambda_{\max}%
(\bar{X})-\varepsilon)\geq0}}(\mu_{1}\lambda_{1}(\bar{X}),\cdots,\mu
_{n}\lambda_{n}(\bar{X}))^{T}%
\end{align*}

\end{dem}

\section{A Technical Appendix: Auxiliary Results and Proofs}

In this appendix we report the auxiliary results which were needed for the
proof of Theorem \ref{tmain}. Recall that $(X,Y)$ denotes a dual pair of
(real) vector spaces with an associated separating bilinear form (dual
pairing)\ denoted by $\left\langle \cdot,\cdot\right\rangle ,\ $so that $X$
and $Y$ are endowed with compatible topologies.

In what follows we fix a lsc convex proper function
\[
\Phi:X\rightarrow\mathbb{R\cup\{+\infty\}},
\]
together with elements
\[
\bar{x}\in\operatorname*{dom}\Phi,\text{ }\delta\geq0,\text{ and }\lambda
\in\left]  -\infty,+\infty\right]  .
\]

The following lemmas give\ estimations for $\mathrm{N}_{[\Phi\leq\lambda
]}^{\delta}(\bar{x})$ under different conditions. The first one uses Slater's
condition at $\lambda$, which means that for some $x_{0}\in X$ we have
$\Phi(x_{0})<\lambda.$

\begin{lem}
\label{lem2}If $\Phi(\bar{x})\leq\lambda<+\infty$ and Slater's condition
holds\ at $\lambda,$ then
\[
\mathrm{N}_{[\Phi\leq\lambda]}^{\delta}(\bar{x})\subset\bigcup_{\mu\geq
0}\partial_{\delta+\mu(\Phi(\bar{x})-\lambda)}\left(  \mu\Phi\right)  (\bar
{x}).
\]

\end{lem}

\begin{dem}
Given\ $\xi\in\mathrm{N}_{[\Phi\leq\lambda]}^{\delta}(\bar{x}),$ we define the
proper lsc convex function $\varphi:X\rightarrow\mathbb{R\cup\{+\infty\}}$ as
\[
\varphi(x):=\max\left\{  \Phi(x)-\lambda,\delta-\left\langle \xi,x-\bar
{x}\right\rangle \right\}  .
\]
From the definition of the $\delta$-normal set we have that $\varphi(x)\geq0$
for all $x\in X.$ Thus, since $\varphi(\bar{x})=\delta,$ it follows that
$\bar{x}$ is a $\delta$-minimum of $\varphi$ and we get, according to
Corollary \ref{finitecase},
\[
\theta\in\partial_{\delta}\varphi(\bar{x})=\bigcup_{\substack{\alpha\in
\lbrack0,1]\\0\leq\eta\leq\alpha\delta+(1-\alpha)(\Phi(\bar{x})-\lambda
)}}\partial_{\eta}((1-\alpha)(\Phi-\lambda)+\alpha(\delta-\xi+\left\langle
\xi,\bar{x}\right\rangle ))(\bar{x}).
\]
In other words, there exist $\alpha\in\lbrack0,1]$ and $\eta\in\lbrack
0,\alpha\delta+(1-\alpha)(\Phi(\bar{x})-\lambda)]$ such that $\alpha\xi
\in\partial_{\eta}(1-\alpha)\Phi(\bar{x}).$ If $\alpha=0,$ then $\eta=0$ (as
$\Phi(\bar{x})\leq\lambda$) and we get $\theta\in\partial\Phi(\bar{x}),$ which
contradicts Slater's condition. So, $\alpha>0$ and the number\ $\mu
:=\frac{1-\alpha}{\alpha}$ $(\geq0)$ is well-defined and satisfies $\frac
{\eta}{\alpha}\leq\delta+\mu(\Phi(\bar{x})-\lambda),$ together with
\[
\xi\in\frac{1}{\alpha}\partial_{\eta}(1-\alpha)\Phi(\bar{x})\subset
\partial_{\frac{\eta}{\alpha}}((1-\alpha)\alpha^{-1})\Phi(\bar{x}%
)\subset\partial_{\delta+\mu(\Phi(\bar{x})-\lambda)}\left(  \mu\Phi\right)
(\bar{x}).
\]

\end{dem}

\bigskip

The following lemma deals with\ the case $\lambda=+\infty.$

\begin{lem}
\label{lem2b}We have that
\[
\mathrm{N}_{\operatorname*{dom}\Phi\cup(\bar{x}+[\Phi\leq\Phi(\bar
{x})]_{\infty})}^{\delta}(\bar{x})=\mathrm{N}_{\operatorname*{dom}\Phi
}^{\delta}(\bar{x})\subset\limsup_{\mu\varepsilon\rightarrow\delta,\text{
}\varepsilon\geq0,\text{ }\mu\downarrow0}\mu\partial_{\varepsilon}\Phi(\bar
{x}).
\]

\end{lem}

\begin{dem}
The first equality is clear because $\operatorname*{dom}\Phi\cup(\bar{x}%
+[\Phi\leq\Phi(\bar{x})]_{\infty})=\operatorname*{dom}\Phi.$ To verify the
inclusion we take $\xi\in\mathrm{N}_{\operatorname*{dom}\Phi}^{\delta}(\bar
{x}).$ We choose $\theta$-neighborhood $V\subset Y$ and $n_{0}\geq1$ such that
$\bar{x},$ $2\bar{x}\in V^{\circ}$ and
\[
\Phi(y)-\Phi(\bar{x})\geq-n\text{ \ for all }y\in V^{\circ}\text{ and all
}n\geq n_{0};
\]
the existence of such a $V$ is due to\ the continuity of the dual pairing,
while $n_{0}$ comes from the weak lower semicontinuity of $\Phi$ and the weak
compactness of $V^{\circ}.$ Fix $n\geq n_{0}.$ Hence, for all $y\in V^{\circ
}\cap\operatorname*{dom}\Phi$ and $\mu\in(0,\frac{\delta}{n^{2}})$ we get
\[
\left\langle \xi,y-\bar{x}\right\rangle \leq\delta\leq\mu(\Phi(y)-\Phi(\bar
{x})+n)+\delta.
\]
The last\ inequality being\ also true when $y\in V^{\circ}\setminus
\operatorname*{dom}\Phi,$ by using the sum rule for the approximate
subdifferential\ (see, e.g., \cite{JBHUPH}) we obtain\ that
\[
\xi\in\partial_{\delta+\mu n}\left(  \mu\Phi+\mathrm{I}_{V^{\circ}}\right)
(\bar{x})\subset\partial_{\delta+2\mu n}(\mu\Phi)(\bar{x})+\partial
_{\delta+2\mu n}\mathrm{I}_{V^{\circ}}(\bar{x})+V.
\]
Since $\bar{x},$ $2\bar{x}\in V^{\circ}$ we easily verify that $\partial
_{\delta+2\mu n}\mathrm{I}_{V^{\circ}}(\bar{x})\subset2(\delta+2\mu
n)V\subset2\delta(1+\frac{2}{n})V,$ and we get\
\[
\xi\in\partial_{\delta(1+\frac{2}{n})}(\mu\Phi)(\bar{x})+2\delta(1+\frac{2}%
{n})V+V=\mu\partial_{\delta\mu^{-1}(1+\frac{2}{n})}\Phi(\bar{x})+2\delta
(1+\frac{2}{n})V+V,
\]
which proves that $\xi\in\limsup_{\mu\varepsilon\rightarrow\delta,\text{
}\varepsilon\geq0,\text{ }\mu\downarrow0}\mu\partial_{\varepsilon}\Phi(\bar
{x}).$
\end{dem}

\begin{lem}
\label{increasing}Assume that $\Phi(\bar{x})\leq\lambda<+\infty.$ Then
\[
\mathrm{N}_{[\Phi\leq\lambda]}^{\delta}(\bar{x})=\bigcap_{\alpha>\delta
}\operatorname*{cl}\left(  \bigcup\nolimits_{\gamma>0}\mathrm{N}_{[\Phi
\leq\lambda+\gamma]}^{\alpha}(\bar{x})\right)  .
\]

\end{lem}

\begin{dem}
We are going to apply Corollary \ref{increasing0} to the (proper, lsc, and
convex) functions $\varphi_{n}:=\mathrm{I}_{[\Phi\leq\lambda+\frac{1}{n}]},$
$n\geq1.$ It is clear that $(\varphi_{n})$ non-decreases as $n$ goes to
$+\infty$ to the function $\varphi:=\mathtt{I}_{[\Phi\leq\lambda]}.$ It is
also clear, since $\varphi_{n}(\bar{x})=\varphi(\bar{x})=0,$ that for each
$\alpha>0$ the sequence of sets $(\partial_{\alpha}\varphi_{n}(\bar{x}))_{n}$
is non-decreasing. Hence,\ Corollary \ref{increasing0} applies and yields
\[
\partial_{\delta}\varphi(\bar{x})=\bigcap_{\alpha>\delta}\limsup
_{n\rightarrow+\infty}\partial_{\alpha}\varphi_{n}(\bar{x})=\bigcap
_{\alpha>\delta}\operatorname*{cl}\left(  \bigcup\nolimits_{n\geq1}%
\partial_{\alpha}\varphi_{n}(\bar{x})\right)  .
\]
Consequently, we write
\begin{align*}
\mathrm{N}_{[\Phi\leq\lambda]}^{\delta}(\bar{x})  &  =\partial_{\delta}%
\varphi(\bar{x})=\bigcap_{\alpha>\delta}\operatorname*{cl}\left(
\bigcup\nolimits_{n\geq1}\partial_{\alpha}\varphi_{n}(\bar{x})\right) \\
&  =\bigcap_{\alpha>\delta}\operatorname*{cl}\left(  \bigcup\nolimits_{n\geq
1}\mathrm{N}_{[\Phi\leq\lambda+\frac{1}{n}]}^{\alpha}(\bar{x})\right)
=\bigcap_{\alpha>\delta}\operatorname*{cl}\left(  \bigcup\nolimits_{\gamma
>0}\mathrm{N}_{[\Phi\leq\lambda+\gamma]}^{\alpha}(\bar{x})\right)  ,
\end{align*}
as we wanted to prove.
\end{dem}

\begin{lem}
\label{lem3}If $\Phi(\bar{x})\leq\lambda<+\infty$, then
\[
\mathrm{N}_{[\Phi\leq\lambda]}^{\delta}(\bar{x})\subset\limsup_{\mu
(\lambda-\Phi(\bar{x})+\varepsilon)\rightarrow\delta,\text{ }\mu>0,\text{
}\varepsilon\geq0}\mu\partial_{\varepsilon}\Phi(\bar{x}).
\]

\end{lem}

\begin{dem}
Take $\xi\in\mathrm{N}_{[\Phi\leq\lambda]}^{\delta}(\bar{x}).$ Suppose first
that Slater's condition holds at $\lambda.$ Then, by\ Lemma \ref{lem2}, there
exists $\mu\geq0$ such that
\[
\xi\in\partial_{\delta+\mu(\Phi(\bar{x})-\lambda)}\left(  \mu\Phi\right)
(\bar{x}).
\]
Hence,$\mathrm{\ }$we are done whenever $\mu>0.$ Otherwise, if $\mu=0,$ then
the last relation reads $\xi\in\partial_{\delta}\left(  0\Phi\right)  (\bar
{x})=\mathrm{N}_{\operatorname*{dom}\Phi}^{\delta}(\bar{x}),$ and Lemma
\ref{lem2b} gives us
\begin{align*}
\xi &  \in\limsup_{\mu\varepsilon\rightarrow\delta,\text{ }\varepsilon
\geq0,\text{ }\mu\downarrow0}\mu\partial_{\varepsilon}\Phi(\bar{x})\\
&  =\limsup_{\mu(\lambda-\Phi(\bar{x})+\varepsilon)\rightarrow\delta,\text{
}\varepsilon\geq0,\text{ }\mu\downarrow0}\mu\partial_{\varepsilon}\Phi(\bar
{x})\subset\limsup_{\mu(\lambda-\Phi(\bar{x})+\varepsilon)\rightarrow
\delta,\text{ }\varepsilon,\text{ }\mu\downarrow0}\mu\partial_{\varepsilon
}\Phi(\bar{x}).
\end{align*}

Suppose now that we don't have Slater's condition at $\lambda,$ so that
$\Phi(\bar{x})=\lambda\leq\Phi(x)$ for all $x\in X,$ and $\bar{x}$ is a
minimum of $\Phi.$ Let us pick a $\theta$-neighborhood $V\subset Y.$ Since
$\mathrm{I}_{[\Phi\leq\Phi(\bar{x})]}=\sup_{n\geq1}\mathrm{I}_{[\Phi\leq
\Phi(\bar{x})+\frac{1}{n}]}$ and the sequence $(\mathrm{I}_{[\Phi\leq\Phi
(\bar{x})+\frac{1}{n}]})_{n}$ is non-decreasing, by\ Lemma \ref{increasing} we
get
\[
\mathrm{N}_{[\Phi\leq\Phi(\bar{x})]}^{\delta}(\bar{x})=\bigcap_{\varepsilon
>0}\operatorname*{cl}\left(  \bigcup\nolimits_{n\in\mathbb{N}}\mathrm{N}%
_{[\Phi\leq\Phi(\bar{x})+\frac{1}{n}]}^{\delta+\varepsilon}(\bar{x})\right)
.
\]
Hence, for every $\varepsilon>0$ small enough there exists $n_{0}\in
\mathbb{N}$ such that $\xi\in\mathrm{N}_{[\Phi\leq\Phi(\bar{x})+\frac{1}{n}%
]}^{\delta+\varepsilon}(\bar{x})+V$ for all $n\geq1.$ As Slater's condition
obviously holds at $\Phi(\bar{x})+\frac{1}{n},$ by\ Lemma \ref{lem2} we obtain
that
\[
\xi\in\bigcup_{\mu\geq0}\partial_{\delta+\varepsilon-\frac{\mu}{n}}\left(
\mu\Phi\right)  (\bar{x})+V\subset\bigcup_{\mu\geq0}\partial_{\delta
+\varepsilon}\left(  \mu\Phi\right)  (\bar{x})+V.
\]
But $\theta\in\partial\Phi(\bar{x}),$ since Slater's condition does not hold
at $\lambda,$ and so\
\[
\xi\in\bigcup_{\mu>0}\partial_{\delta+\varepsilon}\left(  \mu\Phi\right)
(\bar{x})+V=\bigcup_{\mu>0}\mu\partial_{\frac{\delta+\varepsilon}{\mu}}%
\Phi(\bar{x})+V\subset\bigcup_{\mu>0}\mu\partial_{\frac{\delta+\varepsilon
_{0}}{\mu}}\Phi(\bar{x})+V.
\]
The conclusion follows then from the arbitrariness of $\varepsilon>0$ and $V.$
\end{dem}

\bigskip

In Lemma \ref{lem4} below we analyze the case where $\bar{x}$ does not belong
to the set $[\Phi\leq\lambda].$ We shall need the following three lemmas,
which may have their own interest, with the objective to prove that the
operator $S$ given in Lemma \ref{reps1} below is maximal monotone.

\begin{lem}
\label{reps0}Fix $v\in X.$ The function $\varepsilon\rightarrow R(\varepsilon
):=\Phi_{\varepsilon}^{\prime}(\bar{x};v)$\ is non-decreasing and continuous
on $\mathbb{R}_{+}.$
\end{lem}

\begin{dem}
Take\ $0\leq\varepsilon_{0}<\varepsilon_{1}.$ If $f(t):=\Phi(\bar{x}%
+tv)-\Phi(\bar{x}),$ $t>0,$ then $t^{-1}(f(t)+\varepsilon_{0})<t^{-1}%
(f(t)+\varepsilon_{1})$ and, so, $R(\varepsilon_{0})\leq R(\varepsilon_{1});$
that is, $R$ is non-decreasing, and\
\[
\lim_{\varepsilon\downarrow\varepsilon_{0}}R(\varepsilon)=\inf_{\varepsilon
>\varepsilon_{0}}R(\varepsilon)=\inf_{t>0}\inf_{\varepsilon>\varepsilon_{0}%
}\frac{f(t)+\varepsilon}{t}=\inf_{t>0}\frac{f(t)+\varepsilon_{0}}%
{t}=R(\varepsilon_{0}),
\]
showing that $R$ is right-continuous at $\varepsilon_{0}.$ If $\varepsilon
_{0}>0,$ given $n\geq1$ and $\varepsilon<\varepsilon_{0}$ we
choose\ $t_{\varepsilon}>0$ such that $R(\varepsilon)\geq\frac
{f(t_{\varepsilon})+\varepsilon}{t_{\varepsilon}}-\frac{1}{n}.$ Then there is
some $\gamma>0$ such that $t_{\varepsilon}\geq\gamma$ for all $\varepsilon
\in(0,\varepsilon_{0})$ close enough to $\varepsilon_{0}$; such a $\gamma$
exists because for otherwise, the relation $\liminf_{\varepsilon
\uparrow\varepsilon_{0}}t_{\varepsilon}=0$ would yield the following
contradiction (as $f$ is lsc at $0$ and $f(0)=0)$%
\[
+\infty>R(\varepsilon_{0})\geq\liminf_{\varepsilon\uparrow\varepsilon_{0}%
}R(\varepsilon)\geq\liminf_{\varepsilon\uparrow\varepsilon_{0}}\frac
{f(t_{\varepsilon})+\varepsilon}{t_{\varepsilon}}-\frac{1}{n}=+\infty.
\]
Now, writing
\[
R(\varepsilon)\geq\frac{f(t_{\varepsilon})+\varepsilon}{t_{\varepsilon}}%
-\frac{1}{n}=\frac{f(t_{\varepsilon})+\varepsilon_{0}}{t_{\varepsilon}}%
+\frac{\varepsilon-\varepsilon_{0}}{t_{\varepsilon}}-\frac{1}{n}\geq
R(\varepsilon_{0})+\frac{\varepsilon-\varepsilon_{0}}{\gamma}-\frac{1}{n},
\]
we infer\ that $R(\varepsilon_{0})\geq\lim_{\varepsilon\uparrow\varepsilon
_{0}}R(\varepsilon)\geq R(\varepsilon_{0})-\frac{1}{n}.$ Thus, the continuity
of $R$ at $\varepsilon_{0}$ follows as $n$ goes to $+\infty.$
\end{dem}

\bigskip

\begin{lem}
\label{reps1}Given $v\in X,$ we define the set-valued mapping
$S:\mathbb{R\rightrightarrows R}$ as $S(\varepsilon):=\emptyset$ for
$\varepsilon<0,$ and
\[
S(\varepsilon):=\{-\lim_{n\rightarrow+\infty}t_{n}^{-1}\mid\lim_{n\rightarrow
+\infty}t_{n}^{-1}(\Phi(\bar{x}+t_{n}v)-\Phi(\bar{x})+\varepsilon
)=\Phi_{\varepsilon}^{\prime}(\bar{x};v)\},\text{ for }\varepsilon\geq0.
\]
Then the following assertions hold\emph{:}

$(i)$ $S(\varepsilon)\neq\emptyset$ for all $\varepsilon>0;$

$(ii)$ The set $S(0)$ is a possibly empty closed interval of $\mathbb{R}_{-};$

$(iii)$ When\ $S(0)=\emptyset$ there exists $s_{\varepsilon}\in S(\varepsilon
)$ such that $s_{\varepsilon}\rightarrow-\infty$ as $\varepsilon\downarrow0;$

$(iv)$ For every $\varepsilon\geq0$ the set $S(\varepsilon)$ is\ convex and closed.
\end{lem}

\begin{dem}
(i) If $\varepsilon>0,$ then each sequence $(t_{n})$ realizing the infimum in
$R(\varepsilon)$ (recall Lemma \ref{reps0}) must converge (up to a
subsequence) to some $t>0$ (including $t=+\infty$), so that $-\frac{1}{t}\in
S(\varepsilon)$ (with the convention that $\frac{1}{+\infty}=0$).

(ii) If $S(0)$ is a non-empty subset of $\mathbb{R}_{-}$, then it is closed,
by the continuity of function $f=\Phi(\bar{x}+\cdot v)-\Phi(\bar{x}).$ If
$s_{0}$ $(\in S(0)\subset\mathbb{R}_{-})$ is the maximum element in $S(0),$
then there is a sequence $(t_{n})$ of positive numbers such that $s_{0}%
=-\lim_{n\rightarrow+\infty}t_{n}^{-1}$ and $\lim_{n\rightarrow+\infty}%
t_{n}^{-1}(\Phi(\bar{x}+t_{n}v)-\Phi(\bar{x}))=\Phi^{\prime}(\bar{x};v).$ If
$s_{0}<0,$ then for $t_{0}:=\frac{-1}{s_{0}}$ we get
\[
\Phi^{\prime}(\bar{x};v)\leq t^{-1}(\Phi(\bar{x}+tv)-\Phi(\bar{x}))\leq
t_{0}^{-1}(\Phi(\bar{x}+t_{0}v)-\Phi(\bar{x}))=\Phi^{\prime}(\bar{x};v)\text{
\ for all }t\in\left]  0,t_{0}\right]  ;
\]
hence, $\left]  -\infty,s_{0}\right]  \subset S(0).$ If $s_{0}=0,$ then
$t_{n}\rightarrow\infty$ and we obtain
\begin{align*}
\Phi^{\prime}(\bar{x};v)  &  \leq\inf_{t>0}t^{-1}(\Phi(\bar{x}+tv)-\Phi
(\bar{x})\\
&  \leq\sup_{t>0}t^{-1}(\Phi(\bar{x}+tv)-\Phi(\bar{x}))\\
&  =\lim_{n\rightarrow+\infty}t_{n}^{-1}(\Phi(\bar{x}+t_{n}v)-\Phi(\bar
{x}))=\Phi^{\prime}(\bar{x};v);
\end{align*}
that is, $\left]  -\infty,0\right]  \subset S(0).$ Thus, in both cases we have
$\left]  -\infty,s_{0}\right]  \subset S(0).$

(iii) Assume that $S(0)$ is empty. Then, as $\varepsilon\downarrow0,$ there
always exist positive numbers $s_{\varepsilon}\in S(\varepsilon)$ such that
$s_{\varepsilon}\rightarrow-\infty.$ In fact, given an $\varepsilon>0,$ we
choose $t_{\varepsilon}>0$ such that%
\begin{equation}
R(\varepsilon)+\varepsilon=\frac{f(t_{\varepsilon})+\varepsilon}%
{t_{\varepsilon}}\geq\frac{f(t_{\varepsilon})}{t_{\varepsilon}}>R(0);
\label{asr}%
\end{equation}
such an $t_{\varepsilon}$ always exists because, for otherwise, there would
exist $\varepsilon_{i}\downarrow0$ and $t_{i}^{n}\rightarrow_{n}+\infty$ such
that $\lim_{n\rightarrow+\infty}(t_{i}^{n})^{-1}(\Phi(\bar{x}+t_{i}^{n}%
v)-\Phi(\bar{x})+\varepsilon_{i})=\Phi_{\varepsilon_{i}}^{\prime}(\bar{x};v)$
for all $i,$ entailing that $\sup_{t>0}t^{-1}(\Phi(\bar{x}+tv)-\Phi(\bar
{x}))\leq\lim_{n\rightarrow+\infty}(t_{i}^{n})^{-1}(\Phi(\bar{x}+t_{i}%
^{n}v)-\Phi(\bar{x})+\varepsilon_{i})=\Phi_{\varepsilon_{i}}^{\prime}(\bar
{x};v)$ for all $i$. In other words, we get $\sup_{t>0}t^{-1}(\Phi(\bar
{x}+tv)-\Phi(\bar{x}))\leq\Phi^{\prime}(\bar{x};v),$ which gives rise to
$S(0)=\left]  -\infty,0\right]  ,$ a contradiction. Consequently,
(\ref{asr})\ makes sense, so that the vacuity of $S(0)$ together with the
continuity of $f$ lead us to\ $t_{\varepsilon}\rightarrow0^{+}$ (recall Lemma
\ref{reps0}). In other words, $s_{\varepsilon}=-t_{\varepsilon}{}^{-1}$ goes
to $-\infty$ as $\varepsilon$ goes to $0$.

(iv) Since the function $t\rightarrow t^{-1}(f(t)+\varepsilon)$\ $($for
$\varepsilon>0)$ is\ quasi-convex (has convex sublevel sets) and continuous,
the set $A\subset\lbrack0,+\infty]$ defined as
\[
A:=\{t\geq0\mid\exists t_{n}\rightarrow t\text{ s.t. }\lim_{n\rightarrow
+\infty}t_{n}^{-1}(\Phi(\bar{x}+t_{n}v)-\Phi(\bar{x})+\varepsilon
)=R(\varepsilon)\}
\]
is convex and closed. Moreover, $0\notin A$ and the image of $A$ by the
function $\rho(t):=-\frac{1}{t}$ $(t>0)$ coincides with $S(\varepsilon).$
Hence, since function $\rho$ is monotone and continuous we conclude that
$S(\varepsilon)$ is\ convex and closed.
\end{dem}

\bigskip

\begin{lem}
\label{reps}With the notation of Lemmas \emph{\ref{reps0} }and\emph{
\ref{reps1}, }$S$ is a maximal monotone operator.
\end{lem}

\begin{dem}
To show that $S$ is monotone, we pick $(\varepsilon_{i},s_{i})\in S$ (the
graph of $S$), $i=0,1,$ with $0<\varepsilon_{0}<\varepsilon_{1}.$ Then for
each $i=0,1$ there is a sequence $(t_{i}^{n})^{-1}\rightarrow-s_{i}\ $such
that $\lim_{n\rightarrow+\infty}(t_{i}^{n})^{-1}(\Phi(\bar{x}+t_{i}^{n}%
v)-\Phi(\bar{x})+\varepsilon_{i})=R(\varepsilon_{i})$ (recall Lemma
\ref{reps0}); hence $t_{i}^{n}>0.$ Writing
\begin{align*}
R(\varepsilon_{1})  &  =\lim_{n\rightarrow\infty}(t_{\varepsilon_{1}}%
^{n})^{-1}(f(t_{\varepsilon_{1}}^{n})+\varepsilon_{1})\\
&  =\lim_{n\rightarrow\infty}((t_{\varepsilon_{1}}^{n})^{-1}(f(t_{\varepsilon
_{1}}^{n})+\varepsilon_{0})+(t_{\varepsilon_{1}}^{n})^{-1}(\varepsilon
_{1}-\varepsilon_{0}))\\
&  \geq R(\varepsilon_{0})+\liminf_{n\rightarrow\infty}(t_{\varepsilon_{1}%
}^{n})^{-1}(\varepsilon_{1}-\varepsilon_{0})\\
&  =\lim_{n\rightarrow\infty}((t_{\varepsilon_{0}}^{n})^{-1}(f(t_{\varepsilon
_{0}}^{n})+\varepsilon_{0})+\liminf_{n\rightarrow\infty}(t_{\varepsilon_{1}%
}^{n})^{-1}(\varepsilon_{1}-\varepsilon_{0})\\
&  \geq\lim_{n\rightarrow\infty}((t_{\varepsilon_{0}}^{n})^{-1}%
(f(t_{\varepsilon_{0}}^{n})+\varepsilon_{1})+\liminf_{n\rightarrow\infty
}(t_{\varepsilon_{0}}^{n})^{-1}(\varepsilon_{0}-\varepsilon_{1})+\liminf
_{n\rightarrow\infty}(t_{\varepsilon_{1}}^{n})^{-1}(\varepsilon_{1}%
-\varepsilon_{0})\\
&  \geq R(\varepsilon_{1})+\liminf_{n\rightarrow\infty}((t_{\varepsilon_{0}%
}^{n})^{-1}-(t_{\varepsilon_{1}}^{n})^{-1})(\varepsilon_{0}-\varepsilon_{1}),
\end{align*}
we deduce\ that $(\varepsilon_{0}-\varepsilon_{1})(s_{0}-s_{1})\geq0,$ and the
monotonicity of $S$ follows. To check the maximality of $S,$ we observe that
the function $\psi:\mathbb{R}\rightarrow\overline{\mathbb{R}},$ defined as
$\psi(\varepsilon):=\inf\{s\mid s\in S(\varepsilon)\}$ for $\varepsilon\geq0$
and $-\infty$ otherwise, is non-decreasing (and satisfies $\lim_{\varepsilon
\downarrow0}\psi(\varepsilon)=-\infty);$ so, it possesses left and
right-limits $\psi_{-}$ and $\psi_{+}$ everywhere in $\mathbb{R}_{+}.$ Then,
given an $\varepsilon_{0}>0,$ by using \cite[Theorem 2.1.7]%
{zalinesku2000convex} the function $g$ defined on $\mathbb{R}$ as
$g(\tau):=\int_{\varepsilon_{0}}^{\tau}\psi(s)ds$ is a proper lsc convex
function with $\mathbb{R}_{+}\subset\operatorname*{dom}g\subset\overline
{\mathbb{R}}$, and $\partial g(\tau)=\left[  \psi_{-}(\tau),\psi_{+}%
(\tau)\right]  \ $for every $\tau>0,$ while $\partial g(0)=\left]
-\infty,\psi_{+}(0)\right]  $, and $\partial g(\tau)=\emptyset$ for all
$\tau<0.$ Since $S(\varepsilon)$ is convex and closed for every $\varepsilon
\geq0$, by Lemma \ref{reps1}, we infer\ that $\partial g\subset S$ and, so, by
Rockafellar's theorem \cite{rockafellar2015convex} we infer that $S=\partial
g\ $and, in particular, $S$ is maximal monotone.
\end{dem}

\bigskip

Now we are ready to study the set $\mathrm{N}_{[\Phi\leq\lambda]\cup(\bar
{x}+[\Phi\leq\Phi(\bar{x})]_{\infty})}^{\delta}(\bar{x})$ when $\Phi(\bar
{x})>\lambda.$

\begin{lem}
\label{lem4}If $\Phi(\bar{x})>\lambda,$ then
\[
\mathrm{N}_{[\Phi\leq\lambda]\cup(\bar{x}+[\Phi\leq\Phi(\bar{x})]_{\infty}%
)}^{\delta}(\bar{x})\subset\overline{\bigcup_{\mu>0}\partial_{\delta+\mu
(\Phi(\bar{x})-\lambda)}(\mu\Phi)(\bar{x})}.
\]

\end{lem}

\begin{dem}
We pick an element $\xi\in\mathrm{N}_{[\Phi\leq\lambda]\cup(\bar{x}+[\Phi
\leq\Phi(\bar{x})]_{\infty})}^{\delta}(\bar{x})$ such that
\begin{equation}
\xi\notin\overline{\bigcup_{\mu>0}\partial_{\delta+\mu(\Phi(\bar{x})-\lambda
)}(\mu\Phi)(\bar{x})}. \label{ma}%
\end{equation}
Since this last set is convex, nonempty and closed, by Hahn-Banach's Theorem
there exist $v\in X$ and $\alpha\in\mathbb{R}$ such that%
\[
\left\langle \xi,v\right\rangle >\alpha\geq\left\langle x^{\ast}%
,v\right\rangle ,\text{ \ for all}\ x^{\ast}\in\bigcup_{\mu>0}\partial
_{\delta+\mu(\Phi(\bar{x})-\lambda)}(\mu\Phi)(\bar{x});
\]
moreover, because $\partial_{\delta+\mu(\Phi(\bar{x})-\lambda)}(\mu\Phi
)(\bar{x})\supset\partial_{\mu(\Phi(\bar{x})-\lambda)}(\mu\Phi)(\bar{x}%
)=\mu\partial_{\Phi(\bar{x})-\lambda}\Phi(\bar{x})\neq\emptyset,$ we have that
$\alpha\geq0.$ So, one may suppose that $\alpha=\delta,$ so that the
inequalities above read, for all $\mu>0,\ $%
\begin{equation}
\left\langle \xi,v\right\rangle >\delta\geq\left\langle x^{\ast}%
,v\right\rangle ,\text{ \ for all}\ x^{\ast}\in\partial_{\delta+\mu(\Phi
(\bar{x})-\lambda)}(\mu\Phi)(\bar{x})=\mu\partial_{\frac{\delta}{\mu}%
+\Phi(\bar{x})-\lambda}\Phi(\bar{x}); \label{pp}%
\end{equation}
that is $\mu\Phi_{\delta/\mu+\Phi(\bar{x})-\lambda}^{\prime}(\bar{x}%
,v)\leq\delta,$ or, equivalently, for all $\varepsilon\geq0$ (setting
$\varepsilon=\delta/\mu$)
\begin{equation}
\inf_{t>0}\frac{\Phi(\bar{x}+tv)-\lambda+\varepsilon}{t}\leq\varepsilon.
\label{ttp}%
\end{equation}
Let the multifunction $S$ be defined as in Lemma \ref{reps}. First, we assume
that $S(\Phi(\bar{x})-\lambda)\cap\left[  -1,0\right]  \neq\emptyset.$ If
$S(\Phi(\bar{x})-\lambda)$ $(\neq\emptyset\ $by Lemma \ref{reps}) contains a
point $s_{0}\in\left[  -1,0\right[  ,$ then $t_{0}:=\frac{-1}{s_{0}}\geq1$
satisfies $\Phi(\bar{x}+t_{0}v)-\lambda\leq0,$ by (\ref{ttp}), and this leads
us to the following contradiction (recall (\ref{pp})),
\begin{equation}
\left\langle \xi,v\right\rangle =t_{0}^{-1}\left\langle \xi,\bar{x}%
+t_{0}v-\bar{x}\right\rangle \leq t_{0}^{-1}\delta\leq\delta<\left\langle
\xi,v\right\rangle . \label{con}%
\end{equation}
If $S(\Phi(\bar{x})-\lambda)$ contains $0,$ there would exist $t_{n}%
\rightarrow+\infty$ such that\ $\lim_{n\rightarrow+\infty}t_{n}^{-1}(\Phi
(\bar{x}+t_{n}v)-\lambda)=R(\Phi(\bar{x})-\lambda)\leq0,$ which shows that
\[
\sup_{t>0}t^{-1}(\Phi(\bar{x}+tv)-\Phi(\bar{x}))=\lim_{n\rightarrow+\infty
}t_{n}^{-1}(\Phi(\bar{x}+t_{n}v)-\Phi(\bar{x}))=\lim_{n\rightarrow+\infty
}t_{n}^{-1}(\Phi(\bar{x}+t_{n}v)-\lambda)\leq0;
\]
hence, $v\in\lbrack\Phi\leq\Phi(\bar{x})]_{\infty}$ and we get a contradiction
along of (\ref{con}).

Now, we suppose that $S(\Phi(\bar{x})-\lambda)\cap\left[  -1,0\right]
=\emptyset;$ that is, $s<-1$ for all $s\in S(\Phi(\bar{x})-\lambda).$ Then two
cases may occur:

(a) For every $\varepsilon>\Phi(\bar{x})-\lambda$ and $s\in S(\varepsilon)$ we
have $s<-1.$ In this case we pick an $s_{\varepsilon}\in S(\varepsilon)$ and
put $t_{\varepsilon}:=\frac{-1}{s};\ $hence, $t_{\varepsilon}<1,$ so that
\[
\frac{\varepsilon^{-1}(\Phi(\bar{x}+t_{\varepsilon}v)-\lambda)+1}%
{t_{\varepsilon}}\leq1.
\]
Since $\Phi(\bar{x}+\cdot v)$ is bounded from below in $[0,1]$, this last
inequality implies that $t_{\varepsilon}\rightarrow1$ as $\varepsilon
\rightarrow+\infty,$ as well as $\varepsilon^{-1}(\Phi(\bar{x}+t_{\varepsilon
}v)-\lambda)\leq0$ for $\varepsilon$ large enough (because $t_{\varepsilon}%
<1$). Then $\Phi(\bar{x}+v)=\lim_{\varepsilon\rightarrow+\infty}\Phi(\bar
{x}+t_{\varepsilon}v)\leq\lambda$ and we get a contradiction as in (\ref{con}).

(b) There exist some $\varepsilon_{0}>\Phi(\bar{x})-\lambda$ and $s_{0}\in
S(\varepsilon_{0})$\ such that $s_{0}\geq-1.$ Since $S$ is a maximal monotone
operator (Lemma \ref{reps}), it has a convex range and, so, because $s<-1$ for
all $s\in S(\Phi(\bar{x})-\lambda)$ while $s_{0}\geq-1,$ there must exist some
$\varepsilon_{1}>0$ such that\ $-1\in S(\varepsilon_{1});$ that is,
\[
\Phi(\bar{x}+v)-\lambda+\varepsilon_{1}\leq\varepsilon_{1},
\]
and we get $\bar{x}+v\in\lbrack\Phi\leq\lambda]$, which leads us to a
contradiction similar to the one in (\ref{con}). Consequently, (\ref{ma}) is
not true and we must have that $\xi\in\overline{\bigcup_{\mu>0}\partial
_{\delta+\mu(\Phi(\bar{x})-\lambda)}(\mu\Phi)(\bar{x})}.$
\end{dem}

\bigskip

\textbf{Acknowledgment.} We would like to thank the reviewers for
their\ careful reading and for providing valuable suggestions which allowed us
to improve our manuscript.


\begin{thebibliography}{99}                                                                                               %


\bibitem {alvarez2000heavy}\textsc{H. Attouch and F. Alvarez, }\emph{The heavy
ball with friction dynamical system for convex constrained minimization
problems}, in Optimization, Namur (1998), Lecture Notes in Econom. Math.
Systems 481, Springer, Berlin (2000), pp. 25--35.

\bibitem {borwein1982note}\textsc{J. M. Borwein, }\emph{A note on
}$\varepsilon$\emph{-subgradients and maximal monotonicity}, Pacific J. Math.
103(2) (1982), pp. 307--314.

\bibitem {Brondsted72}\textsc{A. Brøndsted}, \emph{On the subdifferential of
the supremum of two convex functions}, Math. Scand. 31 (1972), pp. 225--230.

\bibitem {brondsted1965subdifferentiability}\textsc{A. Brøndsted and R.T.
Rockafellar, }\emph{On the subdifferentiability of convex functions}, Proc.
Amer. Math. Soc. 16 (1965), pp. 605--611.

\bibitem {cabot2009long}\textsc{A. Cabot, H. Engler and S. Gadat}, \emph{On
the long time behavior of second order differential equations with
asymptotically small dissipation}, Trans. Amer. Math. Soc. 361(11) (2009), pp.\ 5983--6017.

\bibitem {cabot2014sequential}\textsc{A. Cabot and L. Thibault},
\emph{Sequential formulae for the normal cone to sublevel sets}, Trans. Amer.
Math. Soc. 366(12) (2014), pp. 6591--6628.

\bibitem {CorreaHantLopez}\textsc{R.\ Correa, A. Hantoute, M. A. López,
}\emph{Towards supremum-sum subdifferential calculus free of qualification
conditions, }SIAM J. Optim., 26(4) (2016), 2219--2234.

\bibitem {HanLopCor16}\textsc{R.\ Correa, A. Hantoute, M. A. López,
}\emph{Weaker conditions for subdifferential calculus of convex functions, }J.
Funct. Anal. 271(5) (2016), pp. 1177--1212.

\bibitem {correa2016bronsted}\textsc{R.\ Correa, A. Hantoute, P. Pérez-Aros,
}\emph{On Brøndsted-Rockafellar Theorem, Maximal Monotonicity of
subdifferential and subdifferential limiting calculus rules for convex lsc
epi-pointed functions in locally convex spaces}, to appear in Math. prog.

\bibitem {ioffe11}\textsc{A. D. Ioffe}, \emph{A note on subdifferentials of
pointwise suprema. }Top 20 (2012), pp. 456-466.

\bibitem {HanLop08}\textsc{A. Hantoute and M. A. López,} \emph{A complete
characterization of the subdifferential set of the supremum of an arbitrary
family of convex functions}. J. Convex Anal. 15 (2008), pp. 831--858.

\bibitem {HanLopZal2008}\textsc{A. Hantoute, M. A. López, and C.
Z\u{a}linescu, }\emph{Subdifferential calculus rules in convex analysis: a
unifying approach via pointwise supremum functions.} SIAM J. Optim. 19 (2008),
pp. 863--882.

\bibitem {hiriart1995subdifferential}\textsc{J-B. Hiriart-Urruty, M.
Moussaoui, A. Seeger, and M. Volle}, \emph{Subdifferential calculus without
qualification conditions, using approximate subdifferentials: a survey},
Nonlinear Anal. 24(12) (1995), pp. 1727--1754.

\bibitem {JBHUPH}\textsc{J.-B. Hiriart-Urruty and R. R. Phelps},
\emph{Subdifferential calculus using }$\varepsilon$\emph{-subdifferentials,
}J. Funct. Anal. 118 (1993), pp. 154-166.

\bibitem {kutateladze}\textsc{S. S. KUTATELADZE}, Convex E-programming, Soviet
Math. Dokl. 20 (1979), pp. 391--393.

\bibitem {lewis}\textsc{A. S. Lewis}, \emph{The convex analysis of unitarily
invariant matrix functions}, J. Convex Anal. 2(1-2) (1995), pp.\ 173--183.

\bibitem {lewisin}\textsc{A. S. Lewis}, \emph{Group invariance and convex
matrix analysis}, SIAM J. Matrix Anal. Appl. 17 (1996), pp. 927--949.

\bibitem {LiNg11}\textsc{C. Li and K. F. Ng}, \emph{Subdifferential Calculus
Rules for Supremum Functions in Convex Analysis,} SIAM J. Optim. 21 (2011),
pp. 782-797.

\bibitem {Volle12}\textsc{M. A. López, M. Volle, }\emph{On the subdifferential
of the supremum of an arbitrary family of extended real-valued functions}.
RACSAM, 105(1) (2011), pp. 3--21.

\bibitem {moreau}\textsc{J. J. Moreau}, Fonctionnelles Convexes. Séminaire sur
les équations aux dérivées partielles. Collège de France 1967.

\bibitem {penot1996subdifferential}\textsc{J. P. Penot,} \emph{Subdifferential
calculus without qualification assumptions}, J. Convex Anal. 2(3) (1996), pp. 207--219.

\bibitem {Rockafellarconjugate}\textsc{R. T. Rockafellar, }Conjugate Duality
and Optimization, Regional Conference Series in Applied Mathematics,
Philadelphia, SIAM, 1974.

\bibitem {rockafellar2009variational}\textsc{R. T. Rockafellar and R. J.B.
Wets}, Variational analysis, Springer-Verlag, Berlin Heideberg, Germany, 1998.

\bibitem {rockafellar2015convex}\textsc{R. T. Rockafellar}, Convex
Analysis\emph{, }Princeton University Press, Princeton, N.J., 1970.

\bibitem {thibault1997sequential}\textsc{L. Thibault, }\emph{Sequential convex
subdifferential calculus and sequential Lagrange multipliers}, SIAM J. Control
Optim. 35(4) (1997), pp. 1434--1444.

\bibitem {zalinesku2000convex}\textsc{C. Z\u{a}linescu}, Convex Analysis in
General Vector Spaces, World Scientific, Singapore, 2002.
\end{thebibliography}
\end{document}